\documentclass[hidelinks,onefignum,onetabnum,final]{required_files/siamart220329}

\usepackage{lipsum}
\usepackage{amsfonts}
\usepackage{graphicx}
\usepackage{algorithmic}
\usepackage{todonotes}
\usepackage{amssymb}
\usepackage{amsmath}
\usepackage{ifthen}
\usepackage{caption}
\usepackage{subcaption}
\usepackage{breakcites}
\usepackage{textgreek}
\usepackage{url}
\usepackage[square,numbers]{natbib}
\usepackage{booktabs}
\usepackage{empheq}

\newcommand*{\colorboxed}{}
\def\colorboxed#1#{%
  \colorboxedAux{#1}%
}
\newcommand*{\colorboxedAux}[3]{%
  \begingroup
    \colorlet{cb@saved}{.}%
    \color#1{#2}%
    \boxed{%
      \color{cb@saved}%
      #3%
    }%
  \endgroup
}
\allowdisplaybreaks
\newcommand{\ra}[1]{\renewcommand{\arraystretch}{#1}}

\newcommand{\displayproofs}{true}
\newcommand{\displayfigures}{true}
\newcommand{\displaycomments}{false}

\ifpdf
  \DeclareGraphicsExtensions{.eps,.pdf,.png,.jpg}
\else
  \DeclareGraphicsExtensions{.eps}
\fi

\newcommand\Rd[2]{\mathbb{R}^{#1\times#2}}

\newcommand\vcol[2]{\begin{pmatrix} #1 \\ #2 \end{pmatrix}}
\newcommand\vrow[2]{\begin{pmatrix} #1 & #2 \end{pmatrix}}
\newcommand\vmat[4]{\begin{pmatrix} #1 & #2 \\ #3 & #4 \end{pmatrix}}
\newcommand\half{\frac12}
\newcommand\pinv[1]{{#1}^\dagger}

\newcommand{\norm}[1]{\left\lVert#1\right\rVert}
\newcommand{\eye}[1]{I_{#1}}
\newcommand{\zeros}[2]{0_{#1, #2}}
\newcommand{\Fmu}[1]{\mathcal{F}_{#1}}

\newcommand{\F}[4]{\mathcal{C}_{#1,#2}^{#3,#4}}

\newcommand{\Fsym}[4]{\mathcal{D}_{#1,#2}^{#3,#4}}
\newcommand{\Fo}[2]{\mathcal{F}_{#1,#2}}

\newcommand{\Co}[2]{\mathcal{C}_{#1}^{#2}}
\newcommand{\Cod}[2]{\tilde{\mathcal{C}}_{#1}^{#2}}
\newcommand{\Cosym}[2]{\mathcal{D}_{#1}^{#2}}
\newcommand{\Cosymd}[2]{\tilde{\mathcal{D}}_{#1}^{#2}}
\newcommand{\prox}[2]{\mathrm{prox}_{#1}\left( #2 \right)}
\newcommand{\proj}[2]{\mathrm{proj}_{#1}\left( #2 \right)}

\newcommand{\Par}[1]{\left( #1 \right)}
\newcommand{\LO}[1]{\mathcal{L}_{#1}}
\newcommand{\LOS}[2]{\mathcal{S}_{#1,#2}}
\newcommand{\LOT}[1]{\mathcal{T}_{#1}}

\newcommand{\proofH}[1]{
    \ifthenelse{\equal{\displayproofs}{false}}
      {
       \textcolor{blue}{\textit{hidden}}
      }
      {
      #1 $~$
      }
}

\newcommand{\figureH}[1]{
    \ifthenelse{\equal{\displayfigures}{false}}
      {
       \textcolor{blue}{\textit{hidden}}
      }
      {
      #1 
      }
}

\newcommand{\commentsH}[1]{
    \ifthenelse{\equal{\displaycomments}{false}}
      {
      }
      {
      #1 
      }
}

\newsiamremark{remark}{Remark}
\newsiamremark{hypothesis}{Hypothesis}
\crefname{hypothesis}{Hypothesis}{Hypotheses}
\newsiamthm{claim}{Claim}
\newsiamremark{property}{Property}
\newsiamthm{conjecture}{Conjecture}

\headers{Interpolation Conditions for Lin. Oper. and applications to PEP}{N. Bousselmi, J. M. Hendrickx and F. Glineur}

\title{Interpolation Conditions for Linear Operators and applications to Performance Estimation Problems\thanks{Submitted to the editors May 26, 2023.
\funding{N. Bousselmi is supported by the French Community of Belgium through a FRIA fellowship (F.R.S-FNRS).}}}

\author{Nizar Bousselmi\thanks{UCLouvain, ICTEAM, Louvain-la-Neuve, Belgium 
  (\email{FirstName.LastName@uclouvain.be}).}
\and Julien M. Hendrickx\footnotemark[2]
\and François Glineur\footnotemark[2]\thanks{UCLouvain, CORE, Louvain-la-Neuve, Belgium }}

\usepackage{amsopn}

\ifpdf
\hypersetup{
  pdftitle={Interpolation Conditions for Linear Operators and applications to Performance Estimation Problems},
  pdfauthor={N. Bousselmi, J. M. Hendrickx and F. Glineur}
}
\fi

\begin{document}

\maketitle

\begin{abstract}
    The Performance Estimation Problem methodology makes it possible to determine the exact worst-case performance of an optimization method. In this work, we generalize this framework to first-order methods involving linear operators. This extension requires an explicit formulation of interpolation conditions for those linear operators. 
    We consider the class of linear operators $\mathcal{M}:x \mapsto Mx$ where matrix $M$ has bounded singular values, and the class of linear operators where $M$ is symmetric and has bounded eigenvalues. 
    We describe interpolation conditions for these classes, i.e. necessary and sufficient conditions that, given a list of pairs $\{(x_i,y_i)\}$, characterize the existence of a linear operator mapping $x_i$ to $y_i$ for all $i$.
    
    Using these conditions, we first identify the exact worst-case behavior of the gradient method applied to the composed objective $h\circ \mathcal{M}$, and observe that it always corresponds to $\mathcal{M}$ being a scaling operator. We then investigate the Chambolle-Pock method applied to $f+g\circ \mathcal{M}$, and improve the existing analysis to obtain a proof of the exact convergence rate of the primal-dual gap. In addition, we study how this method behaves on Lipschitz convex functions, and obtain a numerical convergence rate for the primal accuracy of the last iterate. We also show numerically that averaging iterates is beneficial in this setting.
\end{abstract}

\begin{keywords}
First-order algorithms, Performance estimation, Exact convergence rates, Interpolation conditions, Linear operators, Primal-dual algorithms, Quadratic functions
\end{keywords}

\begin{MSCcodes} 90C25, 90C20, 68Q25, 90C22, 49M29
\end{MSCcodes}

\section{Introduction}

The \textit{Performance Estimation Problem} (PEP) methodology, introduced in \cite{drori2014performance}, computes the exact worst-case performance of a given first-order optimization method on a given class of functions and identifies an instance reaching this worst performance. More precisely, given a method and a performance criterion (lower is better), a PEP is an optimization problem that maximizes this criterion on the application of the method among all possible functions belonging to some class, thus providing the worst possible behavior of the method on the class of functions. It is shown in \cite{taylor2017smooth} that PEPs can be reformulated exactly as semidefinite programs for a wide range of function classes and methods. This provided several tight results on the performance of first-order methods (see e.g. Section \ref{sect:prior_works} for a short review).

In this work, we extend the PEP framework to first-order methods involving linear operators $\mathcal{M} : x \mapsto Mx$, for several classes of matrices $M$ characterized by their eigenvalues or singular values spectrum. The main tool behind this extension is the development of interpolation conditions for linear operators (see Section \ref{sect:interpolation_conditions}). Therefore, we will be able to study the exact worst-case performance of first-order algorithms involving linear operators.

\subsection{Optimization methods involving linear operators} Many methods aim at solving problems involving linear operators in their objective functions or constraints. Iterations of such methods typically themselves call the linear operators of the problem.
We describe below three motivating examples of optimization problems involving linear operators.  

\textbf{Motivating example 1: $\boldsymbol{\min_x g(Mx)}$.} To solve this problem, first-order methods evaluate the gradient of the function $F(x) = g(Mx)$ on some point $x_i$, namely,\vspace{-0.6cm}
\begin{equation}\label{eq:to_decompose}
    \nabla F(x_i) = M^T \nabla g(Mx_i)
\end{equation}

\vspace{-0.2cm}

\noindent where we can see the applications of the matrices $M$ and $M^T$. 

Note that the class of functions $g(Mx)$ includes the quadratic functions $ F(x) = \half x^T Qx$ (when $g(y) = \half ||y||^2$ and $M^TM =Q$). In such a case, the gradient is the linear operator $\nabla F(x_i) = Qx_i$. 
Obtaining the exact worst-case performance of first-order methods on quadratic functions can be done using a known technique (see Section 2.1 in \cite{d2021acceleration} or \cite{ bousselmi2022performance} for recent reviews), which consists in analyzing the roots of a polynomial associated with the method of interest. However, this technique can only deal with a single quadratic objective function, whereas our extension of the PEP methodology can analyze more complex formulations, such as $F(x) = f(x) + \half x^T Q x$ (see \cite{aberdam2022accelerated} for recent work on this class with a nonsmooth function $f$). \smallskip
    
\textbf{Motivating example 2: $\boldsymbol{\min_x f(x) + g(Mx)}$.}  A possible algorithm to solve the problem $\min_x f(x) + g(Mx)$, when $f$ and $g$ are proximable, is the Chambolle-Pock (CP) algorithm \cite{chambolle2011first}. Given step sizes $\tau,\sigma>0$ and a pair of primal-dual points $(x_i,u_i)$, the iteration of CP is
\begin{equation}
    \begin{cases}
        x_{i+1} & = \prox{\tau f}{x_i - \tau M^T u_i}, \\
        u_{i+1} & = \prox{\sigma g^*}{u_i + \sigma M (2x_{i+1}-x_i)}.
    \end{cases}
\end{equation}
The method performs proximal steps, which are already handled in the PEP framework \cite{taylor2018exact}. Again this iteration contains calls to matrices $M$ and $M^T$. \smallskip
    
\textbf{Motivating example 3: $\boldsymbol{\min_{x,y} f(x) + g(y)}$ s.t. $\boldsymbol{M_1 x + M_2 y = b}$.} A famous method to solve the constrained problem $\min_{x,y}  f(x) + g(y)$ s.t. $M_1 x + M_2 y = b$ is the Alternating Direction Method of Multipliers (ADMM) \cite{gabay1976dual} with the following iterations
\begin{equation}
    \begin{cases}
        x_{i+1} \in \arg\min_x f(x) + \frac \rho 2 \norm{M_1 x + M_2 y_i - b + \frac 1 \rho z_i}^2, \\
        y_{i+1} \in \arg\min_y g(y) + \frac \rho 2 \norm{M_1 x_{i+1} + M_2 y - b + \frac 1 \rho z_i}^2, \\
        z_{i+1} = z_i + \rho (M_1 x_{i+1} + M_2 y_{i+1} - b).
\end{cases}
\end{equation}
By contrast to the first two motivating examples, linear operators appear here in the constraint of the problem. However, exactly as in the previous motivating examples, the linear operators end up being used in the iterations of the method. \smallskip

Exploiting additional knowledge about the function structure, e.g. that it can be expressed as $g \circ \mathcal{M}$ or $f+g\circ \mathcal{M} $, improves the analysis of the methods. For an illustrative example, let $g$ be smooth and strongly convex and consider the composed function $g \circ \mathcal{M}$. Despite the fact that $g \circ \mathcal{M}$ only belongs to the class of smooth convex functions (i.e. is not strongly convex), we will see that the gradient method applied to this composed function performs better than on general smooth convex functions (see Section \ref{sect:perf_gMx} and \cite{necoara2019linear}). 
Likewise, CP and ADMM inherently exploit the structure of the problems and cannot be applied to a generic problem $\min_x f(x)$.

More generally, our extension can analyze any first-order method involving linear operators, for example, Primal-Dual Fixed Point \cite{chen2016primal}, Condat-V{\~u} \cite{condat2013primal,vu2013splitting}, Primal-Dual Three-Operator Splitting \cite{yan2018new}, and Proximal Alternating Predictor-Corrector \cite{drori2015simple} algorithms (see \cite{condat2023proximal} for a recent review on these algorithms). These methods solve a wide range of different classical optimization problems, for example, $\ell_2$-regularized robust regression \cite{rousseeuw2005robust}, $\ell_1$-constrained least squares \cite{elden1980perturbation}, basis pursuit \cite{chen1994basis}, total variation deblurring \cite{rudin1992nonlinear,beck2009fast}, and resource allocation \cite{yi2016initialization}.

Existing performance guarantees for the above algorithms are often not tight, may use unusual performance criteria or initial conditions for technical reasons, and are thus difficult to compare. Our extension of PEP remedies these issues by providing the exact worst-case performance of these methods for any of the standard performance criteria and initial conditions.

\subsection{Outline of the paper} In Section \ref{sect:PEP_formulation}, we recall the formulation of the performance estimation problem as a finite-dimensional optimization problem and point out the missing part that we want to extend, i.e.\@ the interpolation conditions for linear operators with bounded eigenvalues or singular values spectrum. In Section \ref{sect:interpolation_conditions}, as the first main contribution, we derive these needed interpolation conditions in an explicit and tractable way to be able to add them to PEP.
As a byproduct, we also present the interpolation conditions for the class of quadratic functions. In Section \ref{sect:exploitation}, as the second main contribution, we exploit our new extension of PEP, first to analyze the gradient method on the problem $\min_x g(Mx)$ and to provide an expression of its worst-case performance. We then tighten the existing performance guarantees on the Chambolle-Pock algorithm and prove the exact convergence rate of the primal-dual gap. We also demonstrate the flexibility of our approach by considering alternate performance criteria and function classes for the same algorithm, for which we obtain numerical performance guarantees.
In Section \ref{sect:conclusion}, we conclude our work and discuss future research directions. \vspace{-0.3cm}

\subsection{Prior PEP work}\label{sect:prior_works}
Since its introduction, the PEP framework has evolved in many directions and settings to analyze a lot of different problems and methods including the gradient method (GM) (and its accelerated version) on smooth strongly or hypoconvex functions \cite{drori2014performance,taylor2017convex,rotaru2022tight}, convex functions with a quadratic upper bound, functions with lower restricted secant inequality and an upper error bound, and non-convex functions satisfying Polyak-Łojasiewicz inequality \cite{abbaszadehpeivasti2023conditions,goujaud2022optimal,guille2022gradient}. More variants of GM have been covered by PEP including GM with exact line search, proximal GM, inexact GM with bounded error on the gradient \cite{de2017worst,taylor2018exact,de2020worst}, and block coordinate descent method \cite{shi2016better, kamri2023worst, abbaszadehpeivasti2022convergence}. Other methods have also been studied by the framework as decentralized methods \cite{colla2023automatic, colla2021automated, colla2022automated}, difference-of-convex-algorithm, the gradient descent-ascent method, the Alternating Direction Method of Multipliers, and Bregman or mirror descent \cite{abbaszadehpeivasti2021rate,abbaszadehpeivasti2022convergence,zamani2023exact,dragomir2021optimal}. In addition to the analysis of optimization problems, the PEP framework can also be used to study problems involving operators like  Halpern iteration for fixed point problems, proximal point algorithm for maximal monotone inclusions, relaxed proximal point algorithm for variational inequalities \cite{lieder2021convergence,gu2020tight,gu2022tight}, splitting methods \cite{ryu2020operator}, and extragradient methods \cite{gorbunov2022extragradient, gorbunov2022last, gorbunov2023convergence}. Some works exploited and extended the PEP framework to develop methodologies that address closely related questions like stochastic optimization, continuous-time model version of first-order methods, and Lyapunov inequalities for first-order methods \cite{taylor2019stochastic,moucer2023systematic,upadhyaya2023automated}. PEP has also been used to optimize and develop optimization methods for different problems as optimal variant of Kelley's cutting-plane method \cite{drori2016optimal}, optimized gradient method for smooth convex functions decreasing the cost function \cite{kim2016optimized,kim2017convergence,kim2018generalizing} and the gradient norm \cite{kim2021optimizing}, optimal gradient method for smooth strongly convex functions \cite{taylor2022optimal}, optimal methods for non-smooth and smooth convex functions \cite{drori2020efficient}, fast iterative shrinkage/thresholding algorithm \cite{kim2018another}, and accelerated proximal points method \cite{kim2021accelerated}. Recent works also proposed to solve non-convex PEP to analyze new problems \cite{gupta2023nonlinear,das2023branch}. Simultaneously to their introduction, some of these extensions have been implemented and are available on the \textit{Matlab} toolbox \texttt{PESTO} \cite{taylor2017performance} and \textit{Python} package \texttt{PEPit} \cite{goujaud2022pepit}. The \href{https://pepit.readthedocs.io/en/latest/index.html}{\texttt{PEPit} documentation website} contains numerous implemented examples of applications of the PEP framework.

Within these works, PEP has been previously used to analyze problems involving linear operators on a few occasions: in \cite{colla2023automatic} for decentralized optimization, in \cite{abbaszadehpeivasti2022convergence} for the random coordinate descent on nonhomogeneous quadratic functions, and in \cite{zamani2023exact} for ADMM. 
We will further detail the contribution of these works at the end of Section \ref{sect:main_results}.
In all cases, they present only potential relaxations of the Performance Estimation Problem for their specific situation, whereas in this work we propose a provably exact formulation of PEP for general problems with linear operators. 

\section{PEP formulation}\label{sect:PEP_formulation}
We summarize the PEP methodology as presented in \cite{taylor2017smooth} and point out the missing parts to analyze problems involving linear operators. As explained, PEP is a framework that analyzes the worst-case behavior of a given optimization method on a given class of functions. For example, a typical PEP could be formulated as follows (but a lot of variations exist). Given the function class $\mathcal{F}$, the optimization method $\mathcal{A}$ performing $N$ iterations, the initial distance $R$, the classical performance criterion $f(x_N) - f(x^*)$ (objective function accuracy after $N$ iterations), and $[N] =\{1,\ldots,N\}$, the PEP is \vspace{-0.25cm}
\begin{equation}\tag{PEP} \label{PEP}
\begin{aligned} 
\max_{x_0,\ldots,x_N,x^*,f} \quad & f(x_N) - f(x^*)\\
\textrm{s.t.} \quad & f \in \mathcal{F}, \\
  & x_i \text{ generated by applying } \mathcal{A} \text{ to } f \text{ from } x_{i-1},~~\forall i\in [N], \\
  & \norm{x_0 - x^*}^2 \leq R^2, \\
  & \norm{\nabla f(x^*)}^2 = 0~\text{(i.e. $x^*$ is optimal)}.
\end{aligned}
\end{equation}
\vspace{-0.35cm} 

\noindent For simplicity, functions in the class $\mathcal{F}$ considered in this example are convex and smooth, so that the optimality condition for $x^*$ is equivalent to stating $\nabla f(x^*)= 0$.

Solving \eqref{PEP} yields the worst-case performance that the method $\mathcal{A}$ can exhibit on a function of the class $\mathcal{F}$ for the performance criterion $f(x_N) - f(x^*)$. Moreover, the maximizer will be an example of worst instance reaching that bound. Note that we can use other performance criteria than $f(x_N) - f(x^*)$, e.g. $||x_N - x^*||^2$, $||\nabla f(x_N)||^2$, $\min_{i\in[N]} ||\nabla f(x_i)||^2$, $f( \frac{1}{N}\sum_{i\in[N]} x_i ) - f(x^*)$, etc.

The conceptual formulation \eqref{PEP} is an infinite-dimensional optimization problem as it involves the class of functions $\mathcal{F}$. However, by discretizing and considering functions $f$ only on the points actually used by the method $\mathcal{A}$, we can rewrite \eqref{PEP} as an equivalent problem with a finite number of variables. Indeed, rather than optimizing over $f\in \mathcal{F}$, we optimize over the points $x_i$, gradients $g_i$ and values $f_i$ that are consistent with a function $f\in \mathcal{F}$ and that can be interpolated by this function.\vspace{-0.1cm}
\begin{definition}[$\mathcal{F}$-interpolability]\label{def:F_interpolability}
Given a function class $\mathcal{F}$, the set of triplets $\{(x_i,g_i,f_i)\}_{i\in [N]}$ is $\mathcal{F}$-interpolable if, and only if,
$
    \exists f\in \mathcal{F}:\begin{cases}
        f(x_i) = f_i, \\
        \nabla f(x_i) = g_i,
    \end{cases} \forall i\in [N].
$
\end{definition}
This definition allows writing the following equivalent discretized formulation of $\eqref{PEP}$\vspace{-0.4cm}
\begin{equation}\tag{PEP-finite} \label{PEP2}
\begin{aligned}
    \max_{S} \quad & f_N - f^*\\
    \textrm{s.t.} \quad & S \text{ is } \mathcal{F}\text{-interpolable}, \\
    & x_i \text{ generated by applying } \mathcal{A} \text{ to } f \text{ from } x_{i-1},~~\forall i\in [N], \\
    & \|x_0 - x^*\|^2 \leq R^2, \\ 
    & \| g^* \|^2 = 0,
\end{aligned}
\end{equation}\vspace{-0.3cm}

\noindent where $S = \{(x_i,g_i,f_i)\}_{i\in \{0\}\cup [N]} \cup \{(x^*,g^*,f^*) \}$. Since we only consider first-order methods, having points $x_i$, gradients $g_i$ and function values $f_i$ is enough, higher-order information on $f$ is not used by the method $\mathcal{A}$.\vspace{-0.1cm}

\subsection{Interpolation conditions}
It remains to express explicitly the first constraint of \eqref{PEP2}. To do so, we need interpolation conditions on the function class $\mathcal{F}$. These are conditions that must satisfy the relevant points $x_i$, $g_i$, and $f_i$ to guarantee that there exists a function $f\in \mathcal{F}$ consistent with those points. In other words, constraints on $\{(x_i,g_i,f_i)\}_{i\in [N]}$ are called interpolation conditions of a class $\mathcal{F}$ when they ensure (and are ensured by) $\{(x_i,g_i,f_i)\}_{i\in [N]}$ being $\mathcal{F}$-interpolable.
For instance, the following theorem provides interpolation conditions for the class $\mathcal{F}_{\mu,L}$ of $L$-smooth $\mu$-strongly convex functions.
\begin{theorem}[\cite{taylor2017smooth}, Theorem 4]
Let $0 \leq \mu < L $ and consider the class $\mathcal{F}_{\mu,L}$. \\
The set of triplets $\{(x_i,g_i,f_i)\}_{i\in [N]}$ is $\mathcal{F}_{\mu,L}$-interpolable if, and only if, $\forall (i,j)\in [N]^2$
\begin{equation}
\label{cond:smooth_convex}
    \small 2\left(1-\frac \mu L \right)\big(f_i - f_j - g_j^T (x_i - x_j)\big) \geq \frac1L\norm{g_i-g_j}^2 + \mu \norm{x_i-x_j}^2- 2\frac\mu L (g_i-g_j)^T(x_i-x_j). 
\end{equation} 
\end{theorem}
It turns out that the nature of interpolation conditions allows in many important cases for a tractable formulation of the PEP, this was shown for instance in \cite{taylor2017smooth} for the class $\Fo{\mu}{L}$. This tractable formulation is a semidefinite problem whose variables are the values $f_i$ and the Gram matrix containing the scalar products between all the iterates $x_i$ and gradients $g_i$. Therefore, the explicit formulation of the interpolation conditions can only involve values $f_i$ and scalar products $x_i^T x_j$, $g_i^T g_j$, and $x_i^Tg_j$ linearly, e.g. \eqref{cond:smooth_convex}, or semidefinite constraints on the Gram matrix of scalar products. Note that since the variables of the problem are the scalar products between the iterates and the gradients, their dimension no longer appears explicitly in the formulation of the problem. We refer the reader to \cite{taylor2017smooth} for more details about the tractable formulation of \eqref{PEP} and interpolation conditions.
 
\subsection{Classes of linear operators}
In this work, we want to extend PEP to methods involving linear operators. As we saw in the motivating examples of the introduction, such first-order methods typically compute the gradient of the objective function but also products of iterates with $M$ and $M^T$. More precisely, we can decompose expression \eqref{eq:to_decompose} of the gradient of a composed function $F(x) = g(Mx)$ as \vspace{-0.2cm}
\begin{align}
    \begin{split}
        y_i & = Mx_i, \\
        v_i =  \nabla F(x_i) \qquad \Leftrightarrow \qquad u_i & = \nabla g(y_i),  \\
        v_i & = M^T u_i .
    \end{split}
\end{align}  \vspace{-0.35cm} 

\noindent Currently, in the PEP framework, it is known how to incorporate equality $u_i = \nabla g(y_i)$. However, $y_i = Mx_i$ and $v_i = M^T u_i$ need new interpolation conditions. Indeed, analyzing such methods with PEP requires the ability to represent the application of a linear operator to a set of points, or even of a linear operator to some points and of the transpose of the same operator to some other points. We formalize this by defining different classes of linear operators $\mathcal{M}$ of interest together with their respective interpolability. In what follows, we will use matrices $M$ and linear operators $\mathcal{M}$ interchangeably as we only work on finite-dimensional spaces. 

We are interested in matrices with bounded eigenvalues or singular values spectrums since the maximal eigenvalue and singular value have an important role in the efficiency of the methods. Such bounds on the spectrum are usually assumed in the literature. We consider general, symmetric, and skew-symmetric matrices, namely,  \vspace{-0.1cm}
\begin{align}
\begin{split}
    \LO{L} & = \{ M: \sigma_{\max}(M) \leq L\}, \\
    \LOS{\mu}{L} & = \{ Q:Q=Q^T,\mu I \preceq Q \preceq LI \}, \\
    \LOT{L} & = \{Q:Q=-Q^T,\sigma_{\max}(Q)\leq L\},
\end{split}
\end{align}  \vspace{-0.2cm}

\noindent with $\sigma_{\max}(M)$ the largest singular value of $M$, and $\preceq$ the L\"owner order. Note that $\LO{L}$ is the class of $L$-Lipschitz linear operators and $\LOS{\mu}{L}$ is the class $\mu$-strongly monotone $L$-Lipschitz symmetric linear operators (see \cite{ryu2020operator} for an analysis on interpolability of general, not necessarily linear operators). We now define operator interpolability in a similar way to function interpolability of Definition \ref{def:F_interpolability}.

\begin{definition}[$\mathcal{L}_{L}$-interpolability]\label{def:R_mat_inter1}
    Sets of pairs $\{(x_i,y_i)\}_{i\in [N_1]}$ and  \linebreak $\{(u_j,v_j)\}_{j\in [N_2]}$ are $\LO{L}$-interpolable if, and only if,\vspace{-0.3cm}
    \begin{equation}
        \exists M\in \LO{L}~:~
        \begin{cases}
            y_i = M x_i, & \forall i\in [N_1], \\
            v_j = M^T u_j, & \forall j\in [N_2].
        \end{cases}
    \end{equation} 
\end{definition}

\begin{definition}[$\LOS{\mu}{L}$-interpolability]\label{def:S_mat_inter1}
    Set of pairs $\{(x_i,y_i)\}_{i\in [N]}$ is $\LOS{\mu}{L}$- interpolable if, and only if,\vspace{-0.3cm}
    \begin{equation}
        \exists Q\in \LOS{\mu}{L}~:~
            y_i = Q x_i,~~ \forall i\in [N]. 
    \end{equation} 
\end{definition}

\begin{definition}[$\LOT{L}$-interpolability]\label{def:Skew_mat_inter1}
    Set of pairs $\{(x_i,y_i)\}_{i\in [N]}$ is $\LOT{L}$- interpolable if, and only if,\vspace{-0.3cm}
    \begin{equation} 
        \exists Q\in \LOT{L}~:~
            y_i = Q x_i,~~ \forall i\in [N]. 
    \end{equation}
\end{definition}

\section{Interpolation conditions for linear operators}\label{sect:interpolation_conditions}
Exploiting PEP to analyze methods involving linear operators requires an explicit formulation of their interpolation conditions. Here we develop tractable necessary and sufficient interpolation conditions for the classes $\LO{L}$, $\LOT{L}$, and $\LOS{\mu}{L}$ of linear operators. We also deduce interpolation conditions for the class of quadratic functions in Section \ref{sect:quad}. 

\subsection{Main results}\label{sect:main_results}
In the sequel, to facilitate manipulation of lists of vectors $\{x_i\}_{i \in [N_1]}$, $\{y_i\}_{i \in [N_1]}$ and $\{u_j\}_{j \in [N_2]}, \{v_j\}_{j \in [N_2]}$, we use notations $X = (x_1~\cdots~x_{N_1})$, $Y = (y_1 ~\cdots ~y_{N_1})$, $U = (u_1 ~\cdots ~u_{N_2})$, $V = (v_1 ~\cdots ~v_{N_2})$ and
write that $(X,Y,U,V)$ is {$\LO{L}$}-interpolable when  $Y=MX$ and $V=M^T U$ for some $M \in \LO{L}$. Similarly, we write that $(X,Y)$ is {$\LOS{\mu}{L}$ or $\LOT{L}$}-interpolable when $Y = QX$ for some $Q$ in $\LOS{\mu}{L}$ or $\LOT{L}$.
We use $\zeros{m}{n}$ for zero matrices (dimension may be omitted) and $||M||=\sigma_{\max}(M)$ for the norm of $M$. We now state our main interpolation theorems.
\begin{theorem}[{$\LO{L}$}-interpolation conditions]\label{th:int_cond_non_sym}
Let $X\in\Rd{n}{N_1}$, $Y\in\Rd{m}{N_1}$, $U\in\Rd{m}{N_2}$, $V\in\Rd{n}{N_2}$, and {$L \geq 0$}.\smallskip

$(X,Y,U,V)$ is {$\LO{L}$}-interpolable if, and only if,
\begin{align} \label{cond:non_symm}
\begin{split}
\begin{cases}
        X^T V = Y^TU, \\
        Y^TY \preceq L^2 X^TX,\\
        V^TV \preceq L^2 U^TU.
\end{cases}
\end{split}
\end{align}

\vspace{-0.25cm} 

\noindent Moreover, if $U=X$ and $V=Y$ (resp. $V=-Y$), then, the interpolant matrix can be chosen symmetric (resp. skew-symmetric).    
\end{theorem}

\begin{corollary}[{$\LOT{L}$}-interpolation conditions]\label{th:int_cond_skewsym}
     Let $X \in \Rd{d}{N}$, $Y\in \Rd{d}{N}$, and {$L \geq 0$}. \smallskip
    
     $(X,Y)$ is {$\LOT{L}$}-interpolable if, and only if,
    \begin{align} \label{cond:skew_symm}
    \begin{split}
    \begin{cases}
        X^T Y = -Y^T X, \\
        Y^TY \preceq L^2 X^TX.
    \end{cases}
    \end{split}
    \end{align}
    
\end{corollary}

\begin{theorem}[{$\LOS{\mu}{L}$}-interpolation conditions]\label{th:int_cond_sym}
    Let $X \in \Rd{d}{N}$, $Y\in \Rd{d}{N}$, and {$-\infty < \mu \leq L < \infty$}. \smallskip
    
    $(X,Y)$ is {$\LOS{\mu}{L}$}-interpolable if, and only if,\vspace{-0.4cm}
    \begin{align} \label{cond:symm}
    \begin{split}
    \begin{cases}
        X^T Y = Y^T X, \\
        (Y-\mu X)^T(LX-Y) \succeq 0.
    \end{cases}
    \end{split}
    \end{align}\vspace{-0.55cm}
      
\end{theorem}
Crucially, these conditions only involve the scalar products between the columns of $(X~U)$ and $(Y~V)$ for {$\LO{L}$} and columns of $(X~Y)$ for {$\LOT{L}$} and {$\LOS{\mu}{L}$}, and they are convex semidefinite-representable constraints on the Gram matrices of these scalar products.

Several additional observations can be made about these conditions. First of all, in Theorem \ref{th:int_cond_non_sym}, the condition $X^TV = Y^T U$ is related to the fact that $X$ and $Y$ are linked by the same matrix (but transposed) as $U$ and $V$. Similarly, in Theorem \ref{th:int_cond_sym} (resp. Corollary \ref{th:int_cond_skewsym}) $X^TY = Y^T X$ (resp. $X^TY = -Y^TX$) is related to the symmetry (resp. skew-symmetry) of the operator. Furthermore, a product $X^TX$ is always symmetric positive semidefinite. Finally, in the nonsymmetric and skew-symmetric cases, $Y^T Y \preceq L^2 X^TX$ is related to the ``maximal amplification" of $X^T X$ through $M^TM$, which produces $Y^TY$, that cannot be larger than $L^2 X^TX$. In the symmetric case, the situation is a bit more sophisticated due to the presence of a nonzero lower bound on the eigenvalues, therefore, we have $Y=QX$ which cannot be ``lower" than $\mu X$ nor ``greater" than $L X$.

Conditions \eqref{cond:non_symm} (without the first equation) and \eqref{cond:symm} were presented as necessary conditions 
for decentralized optimization in Theorem 1 from \cite{colla2023automatic} (in a more specific setting) and in the analysis of ADMM in formulation (13) from \cite{zamani2023exact}. Moreover, in equation (2.7) from \cite{abbaszadehpeivasti2022convergence}, they further restrict the $L$-smooth $\mu$-strongly convex interpolation conditions with an additional necessary condition for nonhomogeneous quadratic functions (see Section \ref{sect:quad} for details). 
Using necessary but not sufficient interpolation conditions in the context of the Performance Estimation Problem still allows obtaining bounds on the worst-case performance, but whose tightness is no longer guaranteed. In this work, we show that conditions \eqref{cond:non_symm} and \eqref{cond:symm} are also sufficient.

\subsection{Proofs of the main results}
We show that the interpolation conditions are indeed necessary and sufficient for $\LO{L}$, $\LOT{L}$, and $\LOS{\mu}{L}$-interpolability. The linear operator that we constructively propose in the nonsymmetric case is actually symmetric (resp. skew-symmetric) when we have the symmetric (resp. skew-symmetric) interpolation conditions. As suggested by the above observations, proving their necessity is much easier than their sufficiency. 

To show that conditions \eqref{cond:non_symm} are sufficient for $\LO{L}$-interpolability, i.e.\@ for the existence of a linear operator $\mathcal{M}$ for $(X,Y,U,V)$, we will first show that under these conditions there exist factorizations $\vrow{X_R}{V_R}$ and $\vrow{Y_R}{U_R}$ of the Gram matrices $G \triangleq \vrow{X}{V}^T \vrow{X}{V}$ and $H \triangleq \vrow{Y}{U}^T \vrow{Y}{U}$, such that $(X_R,Y_R,U_R,V_R)$ is {$\LO{L}$}-interpolable (Lemma \ref{lem:MR_nonsym}). Afterwards, we show how this implies the {$\LO{L}$}-interpolability of the actual vectors $(X,Y,U,V)$ (Lemma \ref{lem:rotation_XYUV}). 
In other words, the proof consists of two steps:
\begin{itemize}
    \item \textbf{Step 1} (Lemma \ref{lem:MR_nonsym}): there exist factorizations $\vrow{X_R}{V_R}$ and $\vrow{Y_R}{U_R}$ of the Gram matrices $G$ and $H$ where $(X_R,Y_R,U_R,V_R)$ is {$\LO{L}$}-interpolable;
    \item \textbf{Step 2} (Lemma \ref{lem:rotation_XYUV}): there exists a rotation mapping lists $(X_R,Y_R,U_R,V_R)$ to the initial vectors $(X,Y,U,V)$.
\end{itemize}

\subsubsection{$\mathcal{L}_{L}$-interpolability of $(X_R,Y_R,U_R,V_R)$ (Step 1)}
Conditions \eqref{cond:non_symm} only involve the scalar products between the columns of $\vrow{X}{V}$ and $\vrow{Y}{U}$, in other words, the Gram matrices $G$ and $H$. Therefore, they apply equally to all sets of vectors leading to the same Gram matrices, i.e. to all factorizations of the Gram matrices.
We first show here that at least one of these factorizations is $\LO{L}$-interpolable.

\begin{lemma}[Existence of $\LO{L}$-interpolable factorizations of Gram matrices]\label{lem:MR_nonsym}
    Let two symmetric positive semidefinite matrices $G= \vmat{A_1}{B_1}{B_1^T}{C_1}$ and $H= \vmat{A_2}{B_2}{B_2^T}{C_2}$ where $A_1, A_2 \in \mathbb{R}^{N_1 \times N_1}$, $C_1, C_2 \in \mathbb{R}^{N_2 \times N_2}$, and $B_1,B_2\in \mathbb{R}^{N_1 \times N_2 }$. \\
    If $G$ and $H$ satisfy\vspace{-0.3cm}
    \begin{equation}
    \begin{cases}
        B_1 = B_2, \\
        A_2 \preceq  L^2 A_1, \\
        C_1 \preceq L^2 C_2,
    \end{cases}
    \end{equation}
    \vspace{-0.42cm}
    
    \noindent then they admit the factorizations \vspace{-0.3cm}
    \begin{align}
        G & = \vrow{X_R}{V_R}^T\vrow{X_R}{V_R}, \\
        H & = \vrow{Y_R}{U_R}^T\vrow{Y_R}{U_R},
    \end{align}
    where \vspace{-0.4cm}
    \begin{align}
        X_R & = \vcol{A_1^\half}{\zeros{N_2}{N_1}} \in \mathbb{R}^{N_1 + N_2 \times N_1}, &  Y_R = \vcol{(\pinv{C_2})^\half B_2^T}{(A_2-B_2 \pinv{C_2}B_2^T)^\half} \in \mathbb{R}^{N_1 + N_2 \times N_1}, \\
        U_R & = \vcol{C_2^\half}{\zeros{N_1}{N_2}} \in \mathbb{R}^{N_2 + N_1 \times N_2}, & V_R  = \vcol{(\pinv{A_1})^\half B_1}{(C_1-B_1^T \pinv{A_1}B_1)^\half} \in \mathbb{R}^{N_2 + N_1 \times N_2},
    \end{align} 
    such that $(X_R,Y_R,U_R,V_R)$ is $\LO{L}$-interpolable.
    \smallskip Moreover, if $A_1=C_2$, $A_2=C_1$, and $B_1 = B_1^T $(resp. $B_1 = -B_1^T$), then, $U_R=X_R$, $V_R=Y_R$ (resp. $V_R=-Y_R)$, and the interpolant matrix is symmetric (resp. skew-symmetric).
    \begin{proof}
    \proofH{The proof consists in explicitly building a matrix $M_R$ defining the operator such that $Y_R = M_R X_R$, $V_R = M_R^T U_R$, and $||M_R||\leq L$. The proof is deferred to Appendix \ref{sect:appendix_proof1}.
    }
    \end{proof}
\end{lemma}
This lemma guarantees that, if $(X,Y,U,V)$ satisfies \eqref{cond:non_symm}, then their associated Gram matrices admit a factorization $(X_R,Y_R,U_R,V_R)$ that is $\LO{L}$-interpolable, but does not yet guarantee the interpolability of the initial vectors $(X,Y,U,V)$ themselves. For example, consider a symmetric case where we just have $(X,Y)$ (hence $N_2=0$) and let  \vspace{-0.1cm}
\begin{equation}
    X = \begin{pmatrix}
    1 \\ 1 \\ 1
\end{pmatrix}, Y = \begin{pmatrix}
    1 \\ 1 \\ 0
\end{pmatrix}\qquad \text{(with $N_1=1$)}
\end{equation}
\vspace{-0.2cm}

\noindent 
which satisfy the $\LOS{\mu}{L}$-interpolation conditions \eqref{cond:symm} for $\mu=0$ and $L=1$. If, given these data $(X,Y)$, one computes the Gram matrix $\begin{pmatrix}X & Y \end{pmatrix}^T \begin{pmatrix}X & Y \end{pmatrix} = \begin{pmatrix} 3 & 2 \\ 2 & 2 \end{pmatrix}$, it turns out that via Lemma \ref{lem:MR_nonsym} one will obtain the following new and different factorization \vspace{-0.2cm}
\begin{equation}
    X_R = \vcol{\sqrt{3}}{0},~~Y_R = \vcol{\frac{2\sqrt{3}}{3}}{\frac{\sqrt{6}}{3}}.
\end{equation}

\vspace{-0.3cm}
\noindent Hence Lemma \ref{lem:MR_nonsym} does not directly guarantee the $\LOS{\mu}{L}$-interpolability of $(X,Y)$, but only of $(X_R,Y_R)$, sharing the same Gram matrix but not necessarily equal to or even having the same dimension as $(X,Y)$. The next lemma guarantees that $(X,Y)$ is also $\LOS{\mu}{L}$-interpolable.

\vspace{-0.5cm}

\subsubsection{Rotation to $(X,Y,U,V)$ (Step 2)}
We now show that if a given factorization of Gram matrices is {$\LO{L}$}-interpolable, then all factorizations of these Gram matrices are {$\LO{L}$}-interpolable.

\begin{lemma}[Rotation between $(X,Y,U,V)$ and $(X_R,Y_R,U_R,V_R)$]\label{lem:rotation_XYUV}
    Let $X_R\in\Rd{n_R}{N_1}$, $Y_R\in\Rd{m_R}{N_1}$, $U_R\in\Rd{m_R}{N_2}$ and $V_R\in\Rd{n_R}{N_2}$.\\

    \vspace{-0.3cm}

    If $(X_R,Y_R,U_R,V_R)$ is {$\LO{L}$}-interpolable, then, $(X,Y,U,V)$ is {$\LO{L}$}-interpolable for all $X\in\Rd{n}{N_1}$, $Y\in\Rd{m}{N_1}$, $U\in\Rd{m}{N_2}$ and $V\in\Rd{n}{N_2}$ such that\vspace{-0.2cm}
    \begin{align}
        \begin{split}
        \vrow{X}{V}^T\vrow{X}{V} & = \vrow{X_R}{V_R}^T\vrow{X_R}{V_R}, \\\vspace{-0.15cm}
        \vrow{Y}{U}^T \vrow{Y}{U} & = \vrow{Y_R}{U_R}^T \vrow{Y_R}{U_R}.
    \end{split}
    \end{align} 
    
    \vspace{-0.2cm}
    
    \noindent Moreover, if $U=X$, $V=Y$ (resp. $V=-Y$), $U_R=X_R$, $V_R=Y_R$ (resp. $V_R=-Y_R$), and interpolant matrix $M_R$ is symmetric (resp. skew-symmetric), then, interpolant matrix $M$ can be chosen symmetric (resp. skew-symmetric).
    \begin{proof}
    \proofH{
    The central idea is that all sets of vectors with the same Gram matrices are equivalent up to a rotation. Some care must however be taken because the dimensions of these vectors are not necessarily the same. The proof is deferred to Appendix \ref{sect:appendix_proof2}.
    }
    \end{proof}
\end{lemma}
We are now able to prove the necessity and sufficiency of interpolation conditions in the nonsymmetric case. The necessity is obtained by a straightforward reasoning and the sufficiency results only on the successive application of Lemmas \ref{lem:MR_nonsym} and \ref{lem:rotation_XYUV}.

\subsubsection{Proof of Theorem \ref{th:int_cond_non_sym}}
\begin{proof}
\proofH{

(Necessity) Let us assume that $(X,Y,U,V)$ is $\LO{L}$-interpolable. First, $Y=M X$ and $V=M^T U$ yield $X^TV = X^T M^T U = Y^T U$. Moreover, $MM^T \preceq L^2 I$ implies $U^T M M^T U \preceq L^2 U^T U$ i.e. $V^TV \preceq L^2 U^TU$ and similarly, $M^TM \preceq L^2 I$ implies $X^T M^T M X \preceq L^2 X^T X$ i.e. $Y^TY \preceq L^2 X^T X$. \smallskip

(Sufficiency) Let us assume that $(X,Y,U,V)$ satisfies conditions \eqref{cond:non_symm}. From Lemma \ref{lem:MR_nonsym}, there exists $(X_R,Y_R,U_R,V_R)$ which is {$\LO{L}$}-interpolable and shares the same Gram matrices as $(X,Y,U,V)$.  Thus, by Lemma \ref{lem:rotation_XYUV}, $(X,Y,U,V)$ is {$\LO{L}$}-inter-polable. Finally, if $U=X$ and $V=Y$ (resp. $V=-Y$), then, we can choose $U_R=X_R$, $V_R=Y_R$ (resp. $V_R=-Y_R$), and $M_R$ symmetric (resp. skew-symmetric) in Lemma \ref{lem:MR_nonsym} and thus $M$ symmetric (resp. skew-symmetric) in Lemma
\ref{lem:rotation_XYUV}. 
}
\end{proof}
The results on the symmetric and skew-symmetric cases come from the nonsymmetric result. In the symmetric case, we apply a shift to also control the lowest eigenvalue.

\subsubsection{Proof of Theorem \ref{th:int_cond_sym}}
\begin{proof}
\proofH{
First, let us define $\Tilde{X}=X$ and $\Tilde{Y} = Y-\frac{L+\mu}{2} X$ and show that requiring $(X,Y)$ to satisfy \eqref{cond:symm} is equivalent to\vspace{-0.2cm}
\begin{align}\label{cond:XYtilde}
\begin{cases}
    \tilde{X}^T \tilde{Y}  = \tilde{Y}^T \tilde{X}, \\
    \tilde{Y}^T \tilde{Y}  \preceq \left(\frac{L-\mu}{2}\right)^2 \tilde{X}^T \tilde{X}.
\end{cases}
\end{align}
\vspace{-0.2cm}

\noindent Indeed, $X^TY = Y^TX \Leftrightarrow X^TY - \frac{L+\mu}{2}X^T X = Y^TX - \frac{L+\mu}{2}X^T X \Leftrightarrow \Tilde{X}^T\Tilde{Y} = \Tilde{Y}^T \Tilde{X}$ and
\vspace{-0.6cm}
\begin{align}
\begin{split}
    (Y-\mu X)^T (Y-LX) \preceq 0 & \Leftrightarrow  \left(\tilde{Y}+ \frac{L-\mu}{2} X\right)^T \left(\tilde{Y} - \frac{L-\mu}{2} X\right) \preceq 0 \\
    & \Leftrightarrow   \tilde{Y}^T \tilde{Y}  \preceq \left(\frac{L-\mu}{2}\right)^2 X^T X.
\end{split}
\end{align}

\vspace{-0.2cm}

\noindent Secondly, from Theorem \ref{th:int_cond_non_sym} in $U=X$ and $V=Y$, conditions \eqref{cond:XYtilde} are equivalent to $(\Tilde{X},\Tilde{Y})$ being {$\LOS{-\frac{L-\mu}{2}}{\frac{L-\mu}{2}}$}-interpolable, therefore, $\Tilde{Y} = \Tilde{Q}\Tilde{X}$ for some $\Tilde{Q} \in \LOS{-\frac{L-\mu}{2}}{\frac{L-\mu}{2}}$. 

Thirdly, $\Tilde{Y} = \Tilde{Q}\Tilde{X}$ with $\Tilde{Q} \in \LOS{-\frac{L-\mu}{2}}{\frac{L-\mu}{2}}$ is equivalent to $Y=QX$ with $Q \in \LOS{\mu}{L}$, indeed, $\tilde{Y}  = \tilde{Q} \tilde{X} \Leftrightarrow Y- \frac{L+\mu}{2}X  = \tilde{Q} X  \Leftrightarrow Y  = (\tilde{Q}+\frac{L+\mu}{2}I) X = QX$.
}
\end{proof}
\subsection{Limiting cases}
Theorems \ref{th:int_cond_non_sym} and \ref{th:int_cond_sym} assume finite values of $L$, we now extend them to ``$L \rightarrow \infty$", that is unbounded singular values and eigenvalues. Again, we are interested in an explicit convex formulation of the conditions. It is not straightforward to take the limit of conditions \eqref{cond:non_symm}, i.e. $\lim_{L\rightarrow \infty}: Y^TY \preceq L^2 X^TX$. Indeed, the constraint must still impose in the limit that the nullspace of $X^TX$ is included in the one of $Y^TY$, which is not directly semidefinite representable. However, it is possible to obtain a tractable formulation of the conditions by considering $\exists L >0 ~:~Y^TY \preceq L^2 X^TX$ instead of the limit. We define $\mathcal{L}$ the class of matrices with arbitrary real singular values and propose the following $\mathcal{L}$-interpolation conditions.
\begin{theorem}[$\mathcal{L}$-interpolation conditions]\label{th:int_cond_non_sym_inf}
Let $X\in\Rd{n}{N_1}$, $Y\in\Rd{m}{N_1}$, $U\in\Rd{m}{N_2}$, and $V\in\Rd{n}{N_2}$.\smallskip

 $(X,Y,U,V)$ is $\mathcal{L}$-interpolable if, and only if,\vspace{-0.2cm}
\begin{align}\label{cond:non_symm_inf}
\begin{split}
\exists L > 0~:~
\begin{cases}
        X^T V = Y^T U, \\
        \vmat{X^TX}{Y^T Y}{Y^T Y}{L^2I} \succeq 0, \\
        \vmat{U^TU}{V^T V}{V^T V}{L^2I} \succeq 0.
\end{cases}
\end{split}
\end{align}
\vspace{-0.2cm}

\begin{proof}
    \proofH{
    By Theorem \ref{th:int_cond_non_sym}, $(X,Y,U,V)$ is $\mathcal{L}$-interpolable if, and only if,\vspace{-0.2cm}
    \begin{equation}
    \begin{cases}
        X^TV & = Y^T U, \\
        \exists L_1>0& :\begin{cases}
                        Y^T Y \preceq L_1^2 X^T X, \\
                        V^T V \preceq L_1^2 U^T U.
                    \end{cases}
    \end{cases}
    \end{equation}  \vspace{-0.3cm}
    Moreover, by Proposition \ref{prop:schur},\vspace{0.1cm}
    \begin{equation}\label{eq:for_note2}
        \exists L_2>0~:~\vmat{X^TX}{Y^T Y}{Y^T Y}{L_2^2I} \succeq 0 \Leftrightarrow  \exists L_2>0~:~(Y^TY)^2 \preceq L_2^2 X^TX.
    \end{equation}
    \vspace{-0.2cm}
    Therefore, we must show that
    \begin{equation}
        \exists L_1>0~:~Y^T Y \preceq L_1^2 X^T X \Leftrightarrow  \exists L_2>0~:~(Y^TY)^2 \preceq L_2^2 X^TX.
    \end{equation}
    Let $C = Y^TY$ and $A = X^T X$\vspace{-0.2cm}
    \begin{align}
    \begin{split}
        \exists L_1>0~:~C \preceq L_1^2 A & \overset{\text{(Prop. \ref{prop:null_loewner})}}{\Leftrightarrow} A\pinv{A}C = C \\
        & \overset{(\text{Prop. \ref{prop:XA_XA})}}{\Leftrightarrow} A\pinv{A}C^2 = C^2  \\
        & \overset{(\text{Prop. \ref{prop:null_loewner}})}{\Leftrightarrow} \exists L_2>0~:~C^2 \preceq L_2^2 A.
    \end{split}
    \end{align}
    And we have the same reasoning for $V^TV \preceq L_1^2 U^T U$.
    }
\end{proof}

\end{theorem}

Existence of $L$ is not an issue for the PEP framework since we can just add a scalar variable to represent $L^2$. Conditions \eqref{cond:non_symm_inf} of Theorem \ref{th:int_cond_non_sym_inf} are thus convex on the Gram matrices $\vrow{X}{V}^T \vrow{X}{V}$, $\vrow{Y}{U}^T \vrow{Y}{U}$ and $L^2$. 

\subsection{Interpolation conditions for quadratic functions} \label{sect:quad}
Let $\mathcal{Q}_{\mu,L}$ the class of homogeneous quadratic functions $f(x) = \half x^T Qx$ where $\mu I \preceq Q \preceq L I$. Our new theorems allow to write the interpolation conditions of $\mathcal{Q}_{\mu,L}$. 

\begin{theorem}[$\mathcal{Q}_{\mu,L}$-interpolation conditions]\label{th:int_cond_quadratic}
Let $-\infty<\mu \leq L < \infty$. \\

The set of triplets $\{(x_i,g_i,f_i)\}_{i\in [N]}$ is $\mathcal{Q}_{\mu,L}$-interpolable if, and only if,\vspace{-0.2cm}
\begin{align} \label{cond:quad_smooth_convex}
    \begin{split}
        \begin{cases}
            X^T G = G^T X, \\
            (G-\mu X)^T (LX-G) \succeq 0, \\
            f_i = \half x_i^T g_i,~~~~\forall i\in [N],
        \end{cases}
    \end{split}
\end{align} \vspace{-0.2cm}
where $X=(x_1 ~\cdots~x_N)$ and $G=(g_1 ~ \cdots ~ g_N)$.
\begin{proof}
    \proofH{
    We have that $\{(x_i,g_i,f_i)\}_{i\in [N]}~\mathcal{Q}_{\mu,L}\text{-interpolable}$ if and only if \vspace{-0.5cm}
    \begin{align*}
         &  \\
    \exists Q\in \mathcal{S}_{\mu,L}~:~ \begin{cases}
        g_i = Q x_i, & \forall i \in [N], \\
        f_i = \half x_i^T Q x_i,& \forall i\in[N],
    \end{cases} & \overset{\text{\phantom{(Th. \ref{th:int_cond_sym})}}}{\Leftrightarrow} 
    \begin{cases}
        g_i = Q x_i, & \forall i \in [N], \\
        f_i = \half x_i^T g_i,& \forall i\in[N],
    \end{cases} \\
    & \overset{\text{(Th. \ref{th:int_cond_sym})}}{\Leftrightarrow}
    \begin{cases}
        X^T G = G^T X \\
        (G-\mu X)^T (LX-G)\succeq 0 \\
        f_i = \half x_i^T g_i, ~~~~\forall i\in[N].
    \end{cases} 
    \end{align*}
    }
\end{proof}
    
\end{theorem}

As we have the inclusion $\mathcal{Q}_{\mu,L} \subseteq \Fo{\mu}{L}$, the quadratic interpolation conditions of Theorem \ref{th:int_cond_quadratic} must imply the general smooth strongly convex interpolation conditions. Indeed, one can also show algebraically that conditions \eqref{cond:quad_smooth_convex} imply conditions \eqref{cond:smooth_convex}.
\begin{lemma} If $\{(x_i,g_i,f_i)\}_{i\in [N]}$ is $\mathcal{Q}_{\mu,L}$-interpolable, then $\{(x_i,g_i,f_i)\}_{i\in [N]}$ is $\Fo{\mu}{L}$-interpolable.
\end{lemma}

As mentioned earlier, the recent work \cite{abbaszadehpeivasti2022convergence} proposed necessary interpolation conditions (originally announced in \cite{drori2019ICCOPT}) for the class of nonhomogeneous quadratic functions $f(x) = \half x^T Qx - b^T x$ (see their Section 2.2). They used the $L$-smooth $\mu$-strongly convex interpolation conditions \eqref{cond:smooth_convex} in addition to the following necessary conditions for nonhomogeneous quadratic functions \vspace{-0.15cm}
\begin{equation} \label{eq:quad_Etienne}
    \half \Par{g_i +g_j}^T \Par{x_i - x_j} = f_i - f_j~~\forall (i,j)\in [N]^2.
\end{equation}

\vspace{-0.15cm}

\noindent These conditions are not sufficient. Indeed, there exist sets of points satisfying both conditions \eqref{cond:smooth_convex} and \eqref{eq:quad_Etienne} and that cannot be interpolated by a quadratic function. For example, the points
\vspace{-0.15cm}
\begin{equation}
    \Par{\vcol{0}{0},\vcol{0}{0},0},~~ \Par{\vcol{1}{0},\vcol{\frac{L+\mu}{2}}{\frac{L-\mu}{2}},\frac{L+\mu}{4}},  \Par{\vcol{-1}{0},\vcol{-\frac{(L+\mu)}{2}}{\frac{L-\mu}{2}},\frac{L+\mu}{4}},
\end{equation}
\vspace{-0.15cm}

\noindent satisfy conditions \eqref{cond:smooth_convex} and \eqref{eq:quad_Etienne}. Yet, there is no function of the form $f(x) = \half x^T Q x - b^T x$ that interpolates these points. Indeed, the first point forces the quadratic to be homogeneous, i.e. $b=0$, and then the linearity of the gradient imposes that the values of the gradient on the two other points are opposite which is impossible for all $\mu\neq L$ (when $\mu=L$, the class of $L$-smooth $L$-strongly convex functions is already the class of nonhomogeneous quadratic functions).

\section{Performance estimation of methods involving linear operators}\label{sect:exploitation}
The new interpolation conditions of Theorems \ref{th:int_cond_non_sym} and \ref{th:int_cond_sym} allow the analysis of any function class currently available in PEP (smooth strongly convex, proximable, hypoconvex, etc.) combined with a linear operator, leading to the exact worst-case performance of methods in these classes.
In practice, we used the \textit{Matlab} toolbox \texttt{PESTO} \cite{taylor2017performance} and added on it our new interpolation conditions. The toolbox is also available on \textit{Python} via the library \texttt{PEPit} \cite{goujaud2022pepit}. Semidefinite programs are solved by the \textit{MOSEK} \cite{mosek} interior-point algorithm.

In this section, we demonstrate the applicability of our extension of the PEP framework. First, we analyze in depth the worst-case performance of the gradient method applied to the first motivating example, $\min_x g(Mx)$ (which also covers problem $\min_x \half x^T Q x$ as a special case). Then, we analyze the more recent and practical Chambolle-Pock algorithm \cite{chambolle2011first}, which was our second motivating example. For the latter, we show the flexibility of our approach with several types of performance guarantees, both in settings found in the literature and in new, previously unstudied settings.

\subsection{Gradient method on $g \circ \mathcal{M}$}
Let us define the class $\F{\mu_g}{L_g}{0}{L_M}$ of functions of the form
\begin{equation}
    F = g \circ \mathcal{M}
\end{equation}
\noindent where $g$ is an $L_g$-smooth $\mu_g$-strongly convex function (where $0<\mu_g \leq L_g$) and $\mathcal{M} : x \mapsto M x$ is defined with a general, not necessarily symmetric matrix $M$ with singular values between $0$ and $L_M$.  
We also define the class $\Fsym{\mu_g}{L_g}{\mu_M}{L_M}$ where $M$ is symmetric with eigenvalues between $\mu_M$ and $L_M$, and $0 \leq \mu_M \leq L_M$.\\

By definition of the classes, we have $\F{\mu_g}{L_g}{0}{L_M} = \F{\mu_g L_M^2}{L_gL_M^2}{0}{1}$ and $\Fsym{\mu_g}{L_g}{\mu_M}{L_M} = \Fsym{\mu_gL_M^2}{L_gL_M^2}{\frac{\mu_M}{L_M}}{1}$, therefore, we will only consider the case $L_M=1$ without loss of generality. For comparison purposes, we will also look at the class $\Fo{\mu_f}{L_f}$ of functions of the form $F = f$ where $f$ is an $L_f$-smooth $\mu_f$-strongly convex function.

We analyze the worst-case performance of the gradient method with fixed step \vspace{-0.25cm}
\begin{equation}\label{eq:GM} \tag{GM}
    x_{k+1} = x_k - \frac{h}{L_F} \nabla F(x_k)   
\end{equation}

\vspace{-0.25cm}

\noindent on the problem $\min_x F(x)$ where $L_F$ is the smoothness constant of the class considered, namely, $L_g$ for $\F{\mu_g}{L_g}{0}{1}$ and $\Fsym{\mu_g}{L_g}{\mu_M}{1}$ and $L_f$ for $\Fo{\mu_f}{L_f}$.

Given a bound $R$ on the initial distance to the solution $||x_0 - x^*||$, we are interested in the worst-case performance $w(\mathcal{F},R,N,\frac{h}{L_F})$ of $N$ iterations of the gradient method with step size $\frac{h}{L_F}$ on the function class $\mathcal{F}$. We define $w(\mathcal{F},R,N,\frac{h}{L_F})$ as the value of the solution of \eqref{PEP} where the method $\mathcal{A}$ is \eqref{eq:GM} with step size $\frac{h}{L_F}$, therefore, it allows to write such following guarantee  \vspace{-0.2cm}
\begin{equation}
     F(x_N) - F(x^*) \leq w\Par{\mathcal{F},R,N,\frac{h}{L_F}} ~~~~\forall F\in \mathcal{F}
 \end{equation}

\vspace{-0.2cm}

\noindent with $x_N$ the iteration $N$ of \eqref{eq:GM} with step size $\frac{h}{L_F}$ on $F$ and $x^*$ the minimizer of $F$. 

The worst-cases $w$ of classes $\Fo{\mu_f}{L_f}$ and $\F{\mu_g}{L_g}{0}{1}$ reduce to simpler cases with the following homogeneity relations, see \cite[Section 4.2.5]{taylor2017convex} for a proof (semicolons are for readability), \vspace{-0.2cm}
\begin{align}\label{eq:homogen}
\begin{split}
    w\left( \Fo{\mu_f}{L_f}; R,N,\frac{h}{L_f} \right) & = L_fR^2 w\left( \Fo{\frac{\mu_f}{L_f}}{1}; 1,N,h \right), \\
    w\left( \F{\mu_g}{L_g}{\mu_M}{1}; R,N,\frac{h}{L_g} \right) & = L_g R^2 w\left( \F{\frac{\mu_g}{L_g}}{1}{\mu_M}{1}; 1,N,h \right).
\end{split}
\end{align}

\vspace{-0.2cm}

\noindent Therefore, without loss of generality, we can consider the cases $L_f = L_g = R = 1$, i.e. $w(\Fo{\mu_f}{1};1,N,h)$ and $w(\F{\mu_g}{1}{\mu_M}{1};1,N,h )$, from which we will deduce the worst-case for the general cases. In the sequel, we will use the following shortened notations (we have the same results and notations for $\Fsym{\mu_g}{L_g}{\mu_M}{1}$) \vspace{-0.1cm}
\[
    w\left(\Fo{\mu_f}{1};1,N,h\right)  = w\left(\Fmu{\mu_f};h\right), \qquad
    w\left(\F{\mu_g}{1}{\mu_M}{1};1,N,h\right)  = w\left(\Co{\mu_g}{\mu_M};h\right).
\]

\vspace{-0.cm}

\noindent Note that $\Fmu{\mu_g}\subseteq \Co{\mu_g}{\mu_M} \subseteq \Fmu{0}$ and $\Fmu{\mu_g}\subseteq \Cosym{\mu_g}{\mu_M} \subseteq \Fmu{\mu_g \mu_M^2}$, therefore, $w(\Fmu{\mu_g};h)\leq w(\Co{\mu_g}{\mu_M};h) \leq w(\Fmu{0};h)$ and $w(\Fmu{\mu_g};h)\leq w(\Cosym{\mu_g}{\mu_M};h) \leq w(\Fmu{\mu_g \mu_M^2};h)$ will always hold. The inclusions merge in $\mu_g=0$.

We are not aware of works addressing the performance of the gradient method on such classes. Reference \cite{necoara2019linear} provided a bound on the performance of gradient method with unit step size on functions $F = g \circ M$ when $M$ has a non-zero lower bound on its singular values, an interesting case but out of our scope of study.

\begin{figure}
     \centering
         \centering
         \figureH{
            \includegraphics[width=0.97\textwidth]{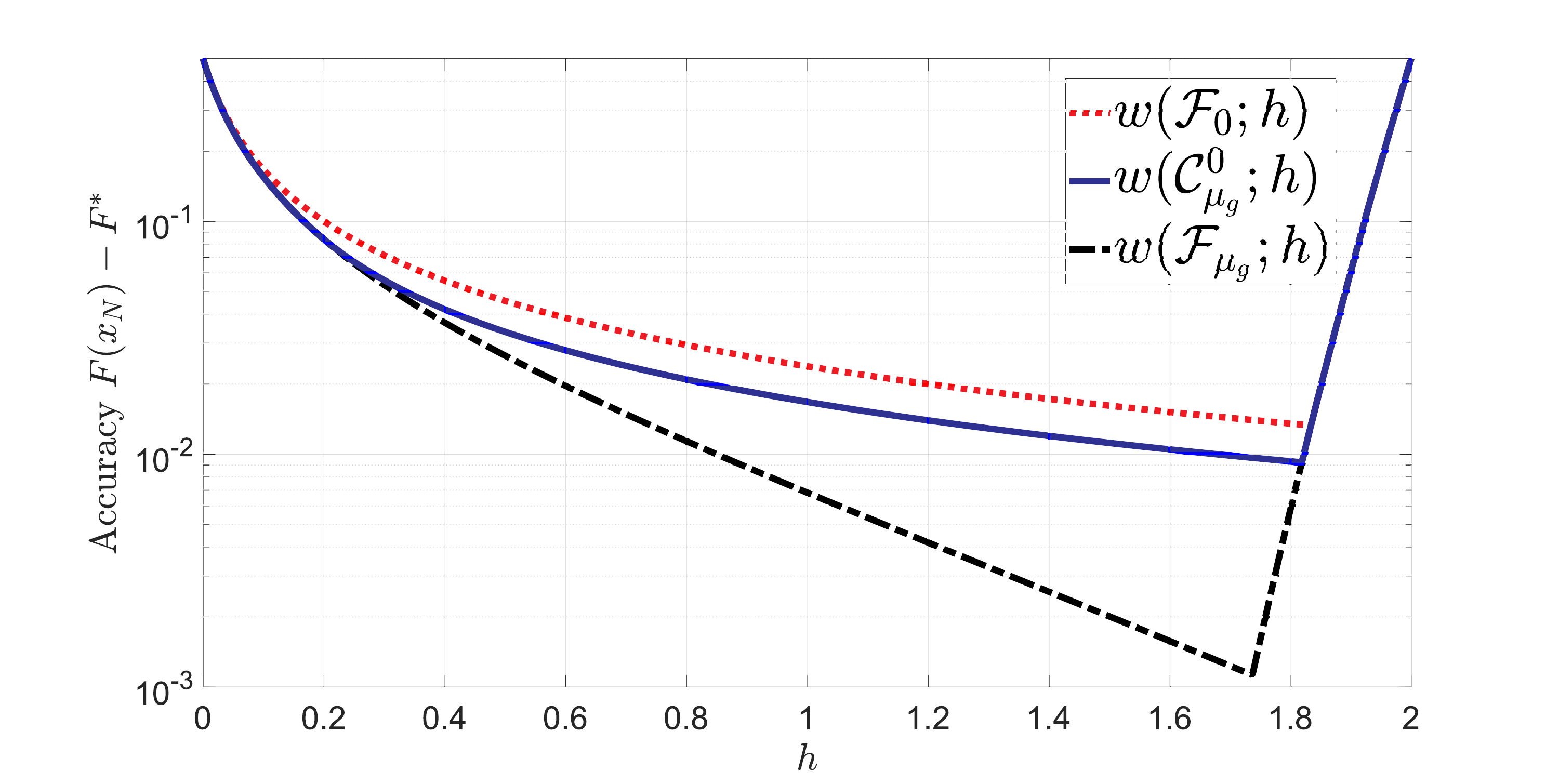}
        }
        \setlength{\belowcaptionskip}{-10pt}
        \caption{Worst-case performance of 10 iterations of \eqref{eq:GM} for varying step size $h\in[0,2]$ on classes $\Fmu{0}$ of $1$-smooth convex functions $f$ (dotted red line), $\Co{0.1}{0}$ of $1$-smooth $0.1$-strongly convex functions $g \circ M$ (solid blue line) and $\Fmu{0.1}$ of $1$-smooth $0.1$-strongly convex functions $f$ (broken black line).}
        \label{fig:exp1and3}
\end{figure}

Solving \eqref{PEP} for the new classes $\Co{\mu_g}{0}$ and $\Cosym{\mu_g}{\mu_M}$ yields  numerical results in Fig. \ref{fig:exp1and3}, namely, the worst-case performance of the gradient method \eqref{eq:GM} after $N=10$ iterations for varying step size $h\in [0,2]$ when applied to classes $\Fmu{0}$, $\Co{0.1}{0}$ and $\Fmu{0.1}$.

After extensive numerical computations with \texttt{PESTO} and analysis of results such that depicted in Fig. \ref{fig:exp1and3}, we were able to identify the exact worst-case performances $w(\Co{\mu_g}{0};h)$ and $w(\Cosym{\mu_g}{\mu_M};h)$. We present the analytical expressions of these worst-cases in Section \ref{sect:perf_gMx}. Before, we propose to review the results of \cite{taylor2017smooth} on the worst-case performance of \eqref{eq:GM} on the class $\Fmu{\mu_f}$ for pedagogical reasons. \\

\subsubsection{Performance on $\Fmu{\mu}$}
Let $N \geq 0$, $h \in [0,2]$, functions $\ell$, $q \in \Fmu{\mu}$ as \vspace{-0.25cm}
\begin{align}\label{eq:def_worst_fun}
\begin{split}
    \ell_{\mu,h}(x) & = 
    \begin{cases}
        \frac{\mu}{2} x^2 + (1-\mu)\tau_{\mu,h}|x| -\left(\frac{1-\mu}{2}\right)\tau_{\mu,h}^2 & \text{if } |x|\geq \tau_{\mu,h}, \\
        \frac{1}{2}x^2 & \text{else},
    \end{cases} \\
    q(x) & = \frac{1}{2} x^2, \\
\end{split}
\end{align}

\vspace{-0.2cm}

\noindent and $\tau_{\mu,h} = \frac{ \mu}{\mu -1 + \Par{1-\mu h }^{-2N}}$. It is conjectured in \cite{taylor2017smooth}, with very strong numerical evidence, that the worst-case performance of \eqref{eq:GM} on $\Fmu{\mu_f}$ is given by \vspace{-0.35cm}
\begin{equation}\label{eq:expr_of_wmuL}
    w\Par{\Fmu{\mu_f};h} = \frac12\max\left\{ \frac{\mu_f}{\mu_f - 1 + (1-\mu_f h)^{-2N}}, (1-h)^{2N} \right\}.
\end{equation}

\vspace{-0.2cm}

\noindent Moreover, as one can check that this worst-case performance corresponds to that of the one-dimensional functions $\ell_{\mu_f,h}$ and $q$, this conjecture also implies that the worst-case $w\Par{\Fmu{\mu_f};h}$ is attained by one-dimensional functions.

\subsubsection{Performance on $\Co{\mu_g}{0}$ and $\Cosym{\mu_g}{\mu_M}$}\label{sect:perf_gMx}
Through our numerous numerical experiments, we also observed that the worst-case functions are one-dimensional. Actually, a more important observation is that, in the worst-case, the operator $\mathcal{M}$ is a scaling (a multiple of the identity), i.e.\@ $M = \alpha I$ for some $\alpha \in \mathbb{R}$. Therefore, we note $\Cod{\mu_g}{0}$ (resp. $\Cosymd{\mu_g}{\mu_M}$) the sub-class of functions $g \circ M$ where $M = \alpha I$ is a scaling operator and propose the following conjecture, supported by our numerical evidence.
\begin{conjecture}\label{conj:scalar}
    Worst-case performances of $\Co{\mu_g}{0}$ and $\Cosym{\mu_g}{\mu_M}$ are reached by scaling operator $M = \alpha I$, i.e., $ w(\Co{\mu_g}{0};h)  = w(\Cod{\mu_g}{0};h)$ and $w(\Cosym{\mu_g}{\mu_M};h) = w(\Cosymd{\mu_g}{\mu_M};h)$.
\end{conjecture}
From now on, given Conjecture \ref{conj:scalar} and the fact that $\Cod{\mu_g}{0} = \Cosymd{\mu_g}{0}$, we will only present the analysis for the symmetric case, as the general case shares the same analysis (but with $\mu_M=0$).
The class of functions $g \circ \alpha I$ can be written as a union of classes of  functions $f$, \vspace{-0.2cm}
\begin{equation}\label{eq:equality_of_sets}
     \Cosymd{\mu_g}{\mu_M} = \bigcup\limits_{\alpha \in [\mu_M,1] } \Fo{\mu_g \alpha^2}{\alpha}^2,
\end{equation}
\vspace{-0.2cm}

\noindent which allows to express the worst-case performance $w(\Cosym{\mu_g}{\mu_M};h)$ thanks to $w(\Fmu{\mu_f};h)$ with
\vspace{-0.2cm}
\newcommand{\myphantom}[0]{\phantom{}}
\begin{align} \label{eq:proof_of_corollary}
\begin{split}
   w\Par{\Cosym{\mu_g}{\mu_M};h} & \overset{\text{Conj. \ref{conj:scalar}}}{=} w \Par{ \Cosymd{\mu_g}{\mu_M};h } \\
    & \overset{\text{\eqref{eq:equality_of_sets}}}{=} w\left( \bigcup\limits_{\alpha\in [\mu_M,1]} \Fo{\mu_g \alpha^2}{ \alpha^2};h \right) \\   
    & \overset{\myphantom}{=} \max_{\alpha\in [\mu_M,1]} w\left( \Fo{\mu_g \alpha^2}{\alpha^2};h \right) \\
         & \overset{\text{Conj. \cite{taylor2017smooth}})}{=} \max_{\alpha\in [\mu_M,1]} w\left( \Fo{\mu_g \alpha^2}{\alpha^2};h \right) \\
    & \overset{\eqref{eq:homogen}}{=} \max_{\alpha\in [\mu_M,1]} \alpha^2 w\left( \Fmu{\mu_g }; \alpha^2 h \right).
\end{split}
\end{align}

\vspace{-0.2cm}

\noindent This reasoning holds only if we know, or conjecture, that the worst-case operator of a method is a scaling of the form $M=\alpha I$. We conjectured it for \eqref{eq:GM} but it is not the case in general. Similarly,  we must not infer that the worst-case operator $M$ is symmetric for all methods and settings. Indeed, one can find an algorithm where the worst-case performance on symmetric operators is strictly better than its performance on general linear operators. 

Given expression \eqref{eq:expr_of_wmuL} of $w( \Fmu{\mu_f};h )$, it is possible to solve the last maximization problem of \eqref{eq:proof_of_corollary} and end up with the following conjecture featuring an explicit convergence rate.

\begin{conjecture}[Worst-case performance $w(\Cosym{\mu_g}{\mu_M};h)$]\label{conj:perf_gMx}
For all $0< \mu_g \leq 1$ and $0\leq \mu_M \leq 1$, we have
\vspace{-0.25cm}
\begin{align}
    w\Par{\Cosym{\mu_g}{\mu_M};h} & = \frac 12 \max\left\{ \frac{\mu_g \alpha^{*^2}}{\mu_g-1+\Par{1-\mu_g \alpha^{*^2} h}^{-2N}}, (1-h)^{2N} \right\}
\end{align}

\vspace{-0.25cm}

\noindent where $\alpha^* = \proj{[\mu_M,1]}{\sqrt{\frac{h_0}{h}}}$ for $h_0$ solution of \vspace{-0.1cm}
\begin{equation}
\begin{cases}
    (1-\mu_g)(1-\mu_g h_0)^{2N+1} = 1 - (2N+1)\mu_g h_0 \\
    0 \leq h_0 \leq \frac{1}{\mu_g}.
\end{cases}
\end{equation}
\end{conjecture}
\vspace{-0.2cm} \noindent
Moreover, we can propose functions reaching this worst-case, therefore, guaranteeing that the worst-case cannot be better (i.e. lower).
\begin{theorem}[Lower bound on the worst-case performance $w(\Cosym{\mu_g}{\mu_M};h)$]\label{th:lower_bound}
    For all $0 < \mu_g \leq 1$ and $0\leq \mu_M \leq 1$, we have \vspace{-0.15cm}
    \begin{align}\label{eq:perf_comp_th}
    w\Par{\Cosym{\mu_g}{\mu_M};h} & \geq \frac 12 \max\left\{ \frac{\mu_g \alpha^{*^2}}{\mu_g-1+\Par{1-\mu_g \alpha^{*^2} h }^{-2N}}, (1-h)^{2N} \right\}
\end{align}

\vspace{-0.25cm}

\noindent with $\alpha^*$ defined in Conjecture \ref{conj:perf_gMx}.
\begin{proof}
    \proofH{
        Functions $\alpha^{*^2}\ell_{\mu_g,\alpha^{*^2}h}(x)$ and $q(x)$, defined in \eqref{eq:def_worst_fun}, belong to $\Cosym{\mu_g}{\mu_M}$ and reach the performance \eqref{eq:perf_comp_th}.
    }
\end{proof}
\end{theorem}
The worst-case performance established in Conjecture \ref{conj:perf_gMx} matches exactly the large number of numerical experiments we performed for many different parameters. Moreover, development \eqref{eq:proof_of_corollary} shows that Conjecture \ref{conj:perf_gMx} relies only on the weaker Conjecture \ref{conj:scalar} and on the previous conjectures of \cite{taylor2017smooth}. We summarize this observation on a corollary. 
\begin{corollary}
    If Conjecture \ref{conj:scalar} holds and $w( \Fmu{\mu_f};h )$ is given by \eqref{eq:expr_of_wmuL} as conjectured in \cite{taylor2017smooth}, then Conjecture \ref{conj:perf_gMx} holds.
\end{corollary}
All these observations and conjectures were made possible thanks to the numerical experiments performed on \texttt{PESTO} with our extension and the extremely helpful insight and information provided by the solution of the different \eqref{PEP} solved.

So far we proposed several partial justifications for Conjecture \ref{conj:perf_gMx} on the worst-case performance of the gradient method. We now give a full proof of the conjecture in a simple case, namely, one iteration of the gradient method with a unit step. \smallskip

\textbf{Proof of Conjecture \ref{conj:perf_gMx} when $N = 1$, $h = 1$, $\mu_M = 0$.}
The optimal primal solution of \eqref{PEP} yields a function reaching the worst-case performance of interest. In contrast, we can extract from the optimal dual solution a combination of inequalities that builds a proof for this performance guarantee (see \cite{goujaud2023fundamental} for more details). Moreover, that type of proof mainly relies on the interpolation conditions. Therefore, solving \eqref{PEP} with our new interpolation conditions will provide a proof that exploits these new conditions. We propose a simple example of such proof for illustration. Analyzing any method with the PEP framework and our interpolation conditions can in principle provide this type of proof.

We consider Conjecture \ref{conj:perf_gMx} in the case $N = 1$, $h = 1$, $\mu_M = 0$, and $\mu_g \leq 0.17$. We still consider the case $L_M=1$ but we do not replace it in the proof for clarity and understanding. The proof below, which shows that  $f(x_1) - f^* \leq \tau ||x_0 - x^*||^2$, was identified from the numerical dual solution of the PEP.
\begin{equation}\hspace{-0.1cm} \label{eq:proof_1_step}
\begin{aligned} 
    & f_1 - f^* \\
    = & 
    \colorboxed{blue}{
    \begin{matrix*}[l]
         \gamma \left(f_1 - f_0 + \nabla g(x_1)^T M(x_0 - x_1) + \frac{\mu_g}{2(1-\mu_g)} ||M(x_0-x_1)||^2\right. \\
         +\left. \frac{1}{2(1-\mu_g)} ||\nabla g(Mx_0)-\nabla g(Mx_1)||^2 - \frac{\mu_g}{1-\mu_g}(\nabla g(Mx_0)-\nabla g(Mx_1))^T M(x_0-x_1)  \right) \\
        + \gamma \left( f_0 - f^*  + \nabla g(Mx_0)^T M(x^* - x_0) + \frac{\mu_g}{2(1-\mu_g)} ||M(x^*-x_0)||^2 \right. \\
        +\left.\frac{1}{2(1-\mu_g)} ||\nabla g(Mx^*)-\nabla g(Mx_0)||^2 -\frac{\mu_g}{1-\mu_g}(\nabla g(Mx^*)-\nabla g(Mx_0))^T M(x^*-x_0) \right) \\
        + (1-\gamma) \left(f_1 - f^* + \nabla g(Mx_1)^T M(x^* - x_1)+ \frac{\mu_g}{2(1-\mu_g)} ||M(x^*-x_1)||^2 \right. \\
        + \left.\frac{1}{2(1-\mu_g)} ||\nabla g(Mx^*)-\nabla g(Mx_1)||^2 - \frac{\mu_g}{1-\mu_g}(\nabla g(Mx^*)-\nabla g(Mx_1))^T M(x^*-x_1)  \right)
    \end{matrix*} }\\
    & \colorboxed{red}{\begin{matrix*}[l]
        + \left( \tau - \frac{1-\gamma}{4} \right) \left( || M(x^* - x_0)||^2 - L_M^2||x^*-x_0||^2\right) \\
        + \frac{1-\gamma}{4} \left(||M (x_0 + x^* - 2 x_1 )||^2 - L_M^2||x^*+x_0 - 2x_1||^2\right) \\
    \end{matrix*}} \\
  & \colorboxed[rgb]{0, 0.6, 0}{- || Pa_{\mu_g} ||^2 - || Pb_{\mu_g} ||^2 - || P c_{\mu_g} ||^2} \\
  & + \tau ||x_0 - x^*||^2 \\
  & \leq \tau ||x_0 - x^*||^2
\end{aligned}
\end{equation}

\noindent where $P = \begin{pmatrix} Mx_0 & Mx_1 & Mx^* & \nabla g(M x_0) & \nabla g(M x_1) & \nabla g(M x^*)\end{pmatrix}$, $\gamma = \frac{1-\mu_g}{2-\mu_g}$, $\tau = \frac{\mu_g}{2(\mu_g-1+(1-\mu_g)^{-2})}$, $x_1 = x_0 - M^T \nabla g(Mx_0)$, and $a_{\mu_g}$, $b_{\mu_g}$, $c_{\mu_g}$ are defined in Appendix \ref{app:proof_1_step}. The expression of $\tau$ corresponds to the one of Conjecture \ref{conj:perf_gMx} for this set of parameters. The equality is an algebraic reformulation that is always satisfied (see Appendix \ref{app:proof_1_step}). The terms in the blue box are nonpositive thanks to the $\mathcal{F}_{\mu_g,1}$-interpolation conditions applied to the function $g$, and the terms in the red box are nonpositive thanks to the $\LO{L_M}$-interpolation conditions applied to the linear operator $M$. Finally, the terms in the green box are a negative sum of squares, thus, always nonpositive. \smallskip

\textbf{Comparison between $\boldsymbol{f}$ and $\boldsymbol{g \circ M}$}. Table \ref{tab:tab1} compares the performance of \eqref{eq:GM} on the classes $\Fmu{\mu_f}$ and $\Cosym{\mu_g}{\mu_M}$ through the performances of $N$ iterations of \eqref{eq:GM} with step size $h$ on the functions $\ell$ and $q$, namely,
\vspace{-0.1cm}
\begin{equation}
\begin{alignedat}{3}\label{eq:p1_p2}
    &p_1(\mu,h)  &&\triangleq \ell_{\mu,h}(x_N) - \ell^*  &&= \frac{1}{2} \frac{\mu}{\mu-1+ \Par{1- \mu h}^{-2N}}, \\
    &p_2(h) &&\triangleq q(x_N)-q^*  &&= \frac{1}{2} \Par{1-h}^{2N}.
\end{alignedat}
\end{equation}

\vspace{-0.8cm}

\begin{table}[H]
    \centering
    \ra{1.5}
    \caption{Worst-case performances and functions of $\Fmu{\mu_f}$ and $\Cosym{\mu_g}{\mu_M}$.}
    \begin{tabular}{@{}lll@{}}
        \toprule
        & $\Fmu{\mu_f}$ & $\Cosym{\mu_g}{\mu_M}$ \\
        \midrule 
        w.c. perf. & $\max\left \{ p_1\left(\mu_f , h\right), p_2(h)\right\}$ & $\max\left \{ \alpha^{*^2} p_1\left(\mu_g, \alpha^{*^2} h\right), p_2(h)\right\}$ \\ [0.3cm] 
        w.c. fun. & $\begin{cases}
            \ell_{\mu_f,h}(x)\\
            q(x)
        \end{cases}$
        & $\begin{cases}
            \alpha^{*^2} \ell_{\mu_g,M^{*^2}h} (x)\\
            q(x) 
        \end{cases}$ \vspace{0.cm} \\ 
        \bottomrule
    \end{tabular}
    \label{tab:tab1}
\end{table}
\noindent There is an interesting difference of performance between the general and composed cases. Let two convex functions $f_1$ and $f_2 = g \circ M$ with $\mu_M = 0$. When the worst-case performance $w$ is given by $p_1$ for $f_2$ and that $h\geq h_0$, then the performance is $\frac{1}{2}\frac{\mu_g \frac{h_0}{h}}{\mu_g-1+(1-\mu_g h_0)^{-2N}} \approx \frac{1}{2} \frac{1}{2Nh+1} e^{-\sqrt{\mu_g}}$ whereas in $p_1$ for $f_1$ it is $\frac{1}{2} \frac{1}{2Nh+1}$. Therefore, there is a gain of around a factor $e^{-\sqrt{\mu_g}}$ between the performance of the gradient method on the convex functions $f_1$ and $f_2$ in this range of values of the step size $h$. \smallskip

\textbf{Optimal step sizes.}
Understanding the exact worst-case performance of \eqref{eq:GM} on $\Cosym{\mu_g}{\mu_M}$ allows to select the optimal step size that minimizes this worst-case performance. Such an optimal design of \eqref{eq:GM} is possible thanks to our extension of PEP. We characterize these optimal step sizes $h \in [0,2]$ minimizing $w( \Cosym{\mu_g}{\mu_M};h )$ from Conjecture \ref{conj:perf_gMx}. 
Optimal steps $h^*(\mu_f)$ of $\Fmu{\mu_f}$ (see \cite{taylor2017smooth}) and $h^*(\mu_g,\mu_M)$ of $\Cosym{\mu_g}{\mu_M}$ satisfy \vspace{-0.2cm}
\begin{align}
    h^*(\mu_f)=h & ~:~ \frac{\mu_f}{\mu_f-1 + \left(1-\mu_f h \right)^{-2N}}  = (1-h)^{2N}, \\
    h^*(\mu_g,\mu_M)=h & ~:~ \frac{\tilde{\mu}(h)}{\mu_g-1 + \left(1-h\tilde{\mu}(h)\right)^{-2N}} = (1-h)^{2N},
\end{align}

\vspace{-0.1cm} 

\noindent where $\tilde{\mu}(h) = \mu_g \proj{[\mu_M^2,1]}{\frac{h_0}{h}}$. Note that both $h^*(\mu_f)$ and $h^*(\mu_g,\mu_M)$ can be easily computed numerically and that they depend on the number of iterations $N$.

Therefore, when facing a $1$-smooth $\mu_f$-strongly convex function $F$, if we do not know anything else about the function, then we should use $h^*\Par{\mu_f}$. However, if we know that the function $F$ can be written as $F = g \circ \mathcal{M}$ where $g$ is $1$-smooth $\mu_g$-strongly convex with $\mu_M \leq ||M|| \leq 1$, then it is preferable to use $h^*\Par{\mu_g,\mu_M}$. \smallskip

\textbf{Performance on quadratic functions.} The class of $L$-smooth, $\mu$-strongly convex, not necessarily homogeneous quadratic functions is easily seen to be equal to  $\mathcal{D}_{1,1}^{\sqrt{\mu},\sqrt{L}}$. Moreover, since the gradient method is invariant to translations, this class shares the same worst-case performance as the class of homogeneous quadratic functions $\mathcal{Q}_{\mu,L}$. Therefore, according to Conjecture \ref{conj:perf_gMx}, the worst-case performance over the class of quadratic functions is\vspace{-0.25cm}
\begin{equation}
   w\Par{ \mathcal{Q}_{\mu,L};h } = w\Par{ \mathcal{D}_{1,1}^{\sqrt{\mu},\sqrt{L}}; h } = \frac{L R^2}{2} \max\left\{q  \left( 1-q h \right)^{2N},(1-h)^{2N} \right\}
\end{equation}

\vspace{-0.25cm}

\noindent where $q = \mathrm{proj}_{[\frac{\mu}{L},1]}( \frac{1}{h(2N+1)})$. Moreover, it is known that the worst-case over the class of quadratic functions must be univariate, a consequence of the polynomial-based analysis mentioned in the introduction (below the first motivating example). Therefore, Conjecture \ref{conj:scalar} holds in the case $L_g=\mu_g$. 

\subsection{Chambolle-Pock method}
We show how to tackle with PEP the analysis of a more sophisticated algorithm, namely the Chambolle-Pock algorithm \cite{chambolle2011first}. 

The Chambolle-Pock algorithm solves problems of the form 
\vspace{-0.15cm}
\begin{equation}\label{eq:prob_CP}
    \min_x f(x) + g(Mx)
\end{equation}

\vspace{-0.25cm}

\noindent where $f$ and $g$ are both convex and proximable and $\mu_M \leq ||M|| \leq L_M$, by applying the following iterations with parameters $\tau>0$ and $\sigma>0$ \vspace{-0.15cm}
\begin{equation}\label{eq:(CP)} \tag{CP}
    \begin{cases}
        x_{i+1} & = \prox{\tau f}{x_i - \tau M^T u_i}, \\
        u_{i+1} & = \prox{\sigma g^*}{u_i + \sigma M (2x_{i+1}-x_i)},
    \end{cases}
\end{equation}

 \vspace{-0.15cm}
 
\noindent where $g^*$ is the convex conjugate function of $g$. 
\subsubsection{Convergence results from the literature}
Existing results from the literature sometimes require very specific assumptions, which makes them difficult to exploit and compare. One of our objectives is to show that the PEP framework allows to obtain guarantees potentially for any set of assumptions, including some that have never been analyzed so far.

For example, the convergence rate stated in the original paper describing the method \cite{chambolle2011first}, reproduced below, requires the existence of sets $B_1$ and $B_2$ "large enough". We note $\mathcal{L}(x,u) = u^T M x + f(x) - g^*(u)$ the Lagrangian of problem \eqref{eq:prob_CP} and $\bar{x}_N = \frac{1}{N}\sum_{i=1}^N x_i$ and $\bar{u}_N = \frac{1}{N}\sum_{i=1}^N u_i$ the averages of the iterates produced by \eqref{eq:(CP)} starting from $x_0$ and $u_0$.
\begin{theorem}[\cite{chambolle2011first}, Theorem 1] \label{th:CP11} Let $f$ and $g$ convex, $||M||\leq L_M$, and $B_1$ and $B_2$ large enough to contain all the iterations $x_i$ and $u_i$ respectively of $\eqref{eq:(CP)}$. If $\tau \sigma L_M^2< 1$, then $\forall N \geq 1$\vspace{-0.2cm}
\begin{equation}
    \mathcal{G}_{B_1 \times B_2} (\bar{x}_N,\bar{u}_N) \leq \frac{D(B_1,B_2)}{N}
\end{equation}

\vspace{-0.2cm}

\noindent where
\vspace{-0.2cm}
\begin{align}
    \begin{split}
        \mathcal{G}_{B_1 \times B_2}(x,u)&  = \max_{u'\in B_2} \mathcal{L}(x,u') - \min_{x'\in B_1} \mathcal{L}(x',u), \\
        D(B_1,B_2) & = \sup_{(x,u)\in B_1 \times B_2} \frac{||x-x_0||^2}{2 \tau} + \frac{||u-u_0||}{2\sigma}.
    \end{split}
    \end{align}
\end{theorem}
The following result, from the same authors, solves this issue and bounds the performance with a quantity that depends only on the initial iterates $x_0$ and $u_0$, but involves some evaluation of the linear operator $M$ of the actual instance of the problem (Remark 2 in \cite{chambolle2016ergodic} adapts the result to remove the dependency in $M$ at the cost of an additional factor).

\begin{theorem}[\cite{chambolle2016ergodic}, Theorem 1]\label{th:CP16}
    Let $f$ and $g$ convex and $||M||\leq L_M$. If $\tau \sigma L_M^2 \leq 1$, then after $N\geq 1$ iterations of the Chambolle-Pock algorithm \eqref{eq:(CP)} started from $x_0$ and $u_0$ we have, for any $x$ and $u$, that the primal-dual gap satisfies
    \begin{equation}\label{eq:th_CP16}
        \small \mathcal{L}(\bar{x}_N,u) - \mathcal{L}(x,\bar{u}_N) \leq \frac{1}{2N} \Par{\frac{\norm{x-x_0}^2}{\tau} + \frac{\norm{u-u_0}^2}{\sigma} - 2 (u-u_0)^TM(x-x_0)}.
    \end{equation}
\end{theorem}
Note that Theorem 2 of \cite{yan2018new} proves a similar result for the Primal-Dual Three-Operator splitting which reduces to the Chambolle-Pock method when the first operator is zero.

Finally, the following result, inspired by a lecture of Prof. Beck (slide 29 of \cite{beck2022}, see Appendix \ref{sect:Beck} for a proof), relies on the previous theorem to bound the primal value accuracy instead of the primal-dual gap. However, the proposed performance guarantee involves a point $\Tilde{u}_N$ which depends on the average iterate $\bar{x}_N$, and which cannot be easily bounded a priori.
\begin{theorem}\label{th:Beck} Let $f$ and $g$ convex and $||M||\leq L_M$. If $\tau \sigma L_M^2 \leq 1$, then after $N\geq 1$ iterations of the Chambolle-Pock algorithm \eqref{eq:(CP)} started from $x_0$ and $u_0$ we have, for any $x$ and $u$, that the primal-dual gap satisfies\vspace{-0.2cm}
    \begin{equation}\label{eq:bound_AB} 
    \small{F(\bar{x}_N) - F(x^*) \leq         \frac{1}{2N} \left(\frac{||x^* - x_0||^2}{\tau} + \frac{|| \tilde{u}_N-u_0||^2}{\sigma} -2(\tilde{u}_N-u_0)^T M (x^*-x_0) \right)}
    \end{equation}
    
    \vspace{-0.15cm} 
    
    \noindent where $F(x) = f(x) + g(Mx)$ and $\tilde{u}_N \in \partial g(M \bar{x}_N)$. 
\end{theorem}

\subsubsection{Convergence results obtained with the new PEP framework}

\begin{figure}[H]
    \centering
    \figureH{
    \includegraphics[width=0.98\textwidth]{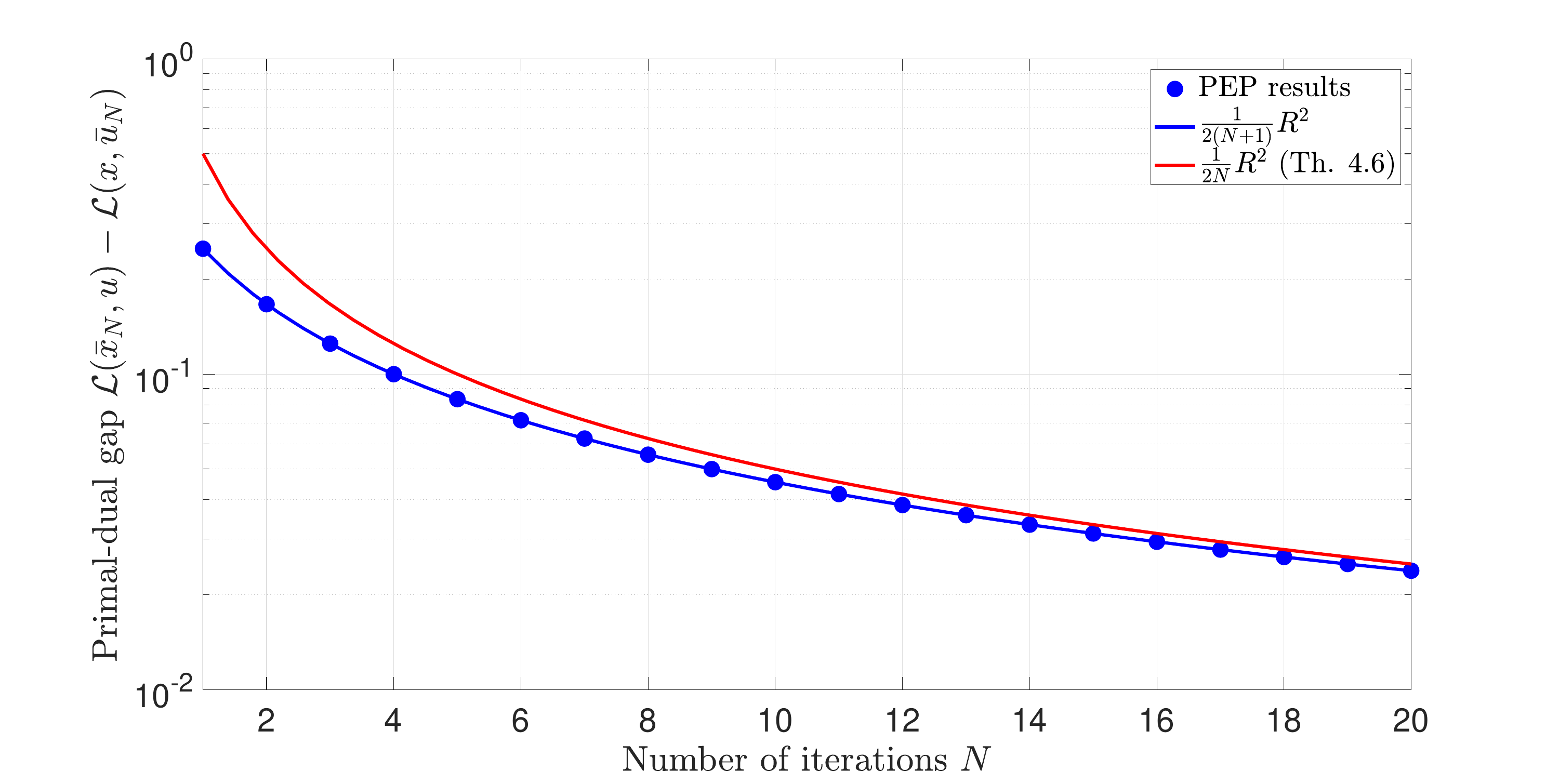}
    }
    \caption{Worst-case performance obtained by our extension of PEP for $N$ iterations of the Chambolle-Pock algorithm \eqref{eq:(CP)} with step size $\tau\sigma L_M^2 \leq 1$ on the problem $\min_x F(x)$ where $F = f + g \circ M$, $f$ and $g$ are convex proximable and $M$ is such that $0\leq ||M||\leq 1$. The performance criterion is the primal-dual gap $\mathcal{L}(\bar{x}_N,u) - \mathcal{L}(x,\bar{u}_N)$ and the initial distance is $R^2 = \frac{\norm{x-x_0}^2}{\tau} + \frac{\norm{u-u_0}^2}{\sigma} - 2 (u-u_0)^TM(x-x_0)$. PEP results (blue dots) are compared to bound \eqref{eq:th_CP16} of Theorem \ref{th:CP16} (red line).
    }
    \label{fig:CP16_improved}
\end{figure}

Theorems \ref{th:CP16} and \ref{th:Beck} can be analyzed using our extension of the PEP framework. First, to study the setting of Theorem \ref{th:CP16}, we impose that the initial distance $R^2 = \frac{\norm{x-x_0}^2}{\tau} + \frac{\norm{u-u_0}^2}{\sigma} - 2 (u-u_0)^TM(x-x_0)$ in the right-hand side of \eqref{eq:th_CP16} is bounded by one, and compute the worst case for the left-hand side (i.e. its maximum), giving us the best possible value of the leading factor. We observe from the tight numerical results provided by PEP that the exact rate is strictly better than the existing bound $\frac{R^2}{2N}$, as shown in Fig. \ref{fig:CP16_improved}. Note that these numerical results were identical for all values of $\tau$ and $\sigma$ such that $\tau\sigma L^2_M \le 1$. We also observe numerically that the method no longer converges when $\tau \sigma L_M^2 > 1$. \vspace{-0.4cm}

Moreover, we were able to identify a simple analytical expression matching the tight numerical bound returned by PEP: it corresponds to $\frac{R^2}{2(N+1)}$, i.e.\@ the same bound as \eqref{eq:th_CP16} but with $N+1$ instead of $N$ (shown as a blue line on Fig. \ref{fig:CP16_improved}).

\noindent Based on the hint, provided by PEP, that this better bound actually holds, we were able to adapt the proof of Theorem \ref{th:CP16} in \cite{chambolle2016ergodic} to obtain the following Theorem, that gives the exact convergence rate.
\begin{theorem} \label{th:CP16_improved}
Let $f$ and $g$ convex and $||M||\leq L_M$. If $\tau \sigma \leq \frac{1}{L_M^2}$, then after $N\geq 1$ iterations of the Chambolle-Pock algorithm \eqref{eq:(CP)} started from $x_0$ and $u_0$ we have, for any $x$ and $u$, that the primal-dual gap satisfies \vspace{-0.cm}
\begin{equation}\label{eq:th_CP16_improved}
    \mathcal{L}(\bar{x}_N,u) - \mathcal{L}(x,\bar{u}_N) \leq \frac{1}{2(N+1)} \Par{\frac{\norm{x-x_0}^2}{\tau} + \frac{\norm{u-u_0}^2}{\sigma} - 2 (u-u_0)^TM(x-x_0)}.
\end{equation}
\vspace{-0.1cm}
\begin{proof}
    Define $||z||^2_A = \frac{||z_x||^2}{\tau} + \frac{||z_u||^2}{\sigma} - 2 z_u^T M z_x$ with $z=(z_x,z_u)$, so that the initial distance becomes $R^2 = \|z-z_0\|_A^2$ with $z=(x,u)$ and $z_0=(x_0,u_0)$. Letting now $z_n = (x_n,u_n)$,
    we sum the following inequality (proved as (15) in \cite{chambolle2011first})
    \[ \mathcal{L}(x_{n+1},u) - \mathcal{L}(x,u_{n+1}) \leq \|z-z_n\|_A^2 - \|z-z_{n+1} \|_A^2 - \| z_{n+1} - z_n \|_A^2 \]
    from $n=0$ to $n=N-1$, without removing the negative terms, yielding \vspace{-0.15cm}
    \begin{align*}
        \sum_{n=0}^{N-1} \mathcal{L}(x_{n+1},u) - \mathcal{L}(x,u_{n+1}) & \leq \frac{||z_0-z||_A^2}{2} - \frac{||z-z_N||_A^2}{2} - \sum_{n=0}^{N-1} \frac{||z_{n+1}-z_n||_A^2}{2} \\
        & \leq \frac{||z_0-z||_A^2}{2} - \frac{1}{N+1}\frac{||z_0-z||_A^2}{2}
    \end{align*}
    where to prove the second inequality we use the following consequence of the convexity of the squared norm, with $s_i=z_i$ for all $0 \le i \le N$ and $s_{N+1} = z$:
            \[ \tfrac{1}{N+1} \|s_0-s_{N+1}\|^2 = (N+1) \|\tfrac{1}{N+1} \textstyle\sum_{n=0}^N (s_{n+1}-s_n)\|^2 \leq  \sum_{n=0}^{N} \|s_{n+1}-s_n\|^2. \]    
        Finally, by convexity of $\mathcal{L}(\cdot,u)-\mathcal{L}(x,\cdot)$, we have \vspace{-0.15cm}
    \begin{align*}
        N (\mathcal{L}(\bar{x}_N,u) - \mathcal{L}(x,\bar{u}_N) ) \leq \sum_{n=0}^{N-1} \mathcal{L}(x_{n+1},u) - \mathcal{L}(x,u_{n+1}) \leq \Par{1-\frac{1}{N+1}} \frac{||z_0-z||_A^2}{2}
    \end{align*}
    leading finally to $\mathcal{L}(\bar{x}_N,u) - \mathcal{L}(x,\bar{u}_N) \leq \frac{1}{N+1} \frac{||z_0-z||_A^2}{2}$.
\end{proof}
\end{theorem}

Instead of bounding the quantity $R^2 = \frac{\norm{x-x_0}^2}{\tau} + \frac{\norm{u-u_0}^2}{\sigma} - 2 (u-u_0)^TM(x-x_0)$, it is also possible with PEP to compute  rates that depend on the simpler squared distance \vspace{-0.15cm}
\begin{equation}\label{eq:new_init}
    R_0^2 = \norm{x-x_0}^2 + \norm{u-u_0}^2.
\end{equation}

\vspace{-0.4cm}

\begin{figure}[H]
    \centering
    \figureH{
    \includegraphics[width=1\textwidth]{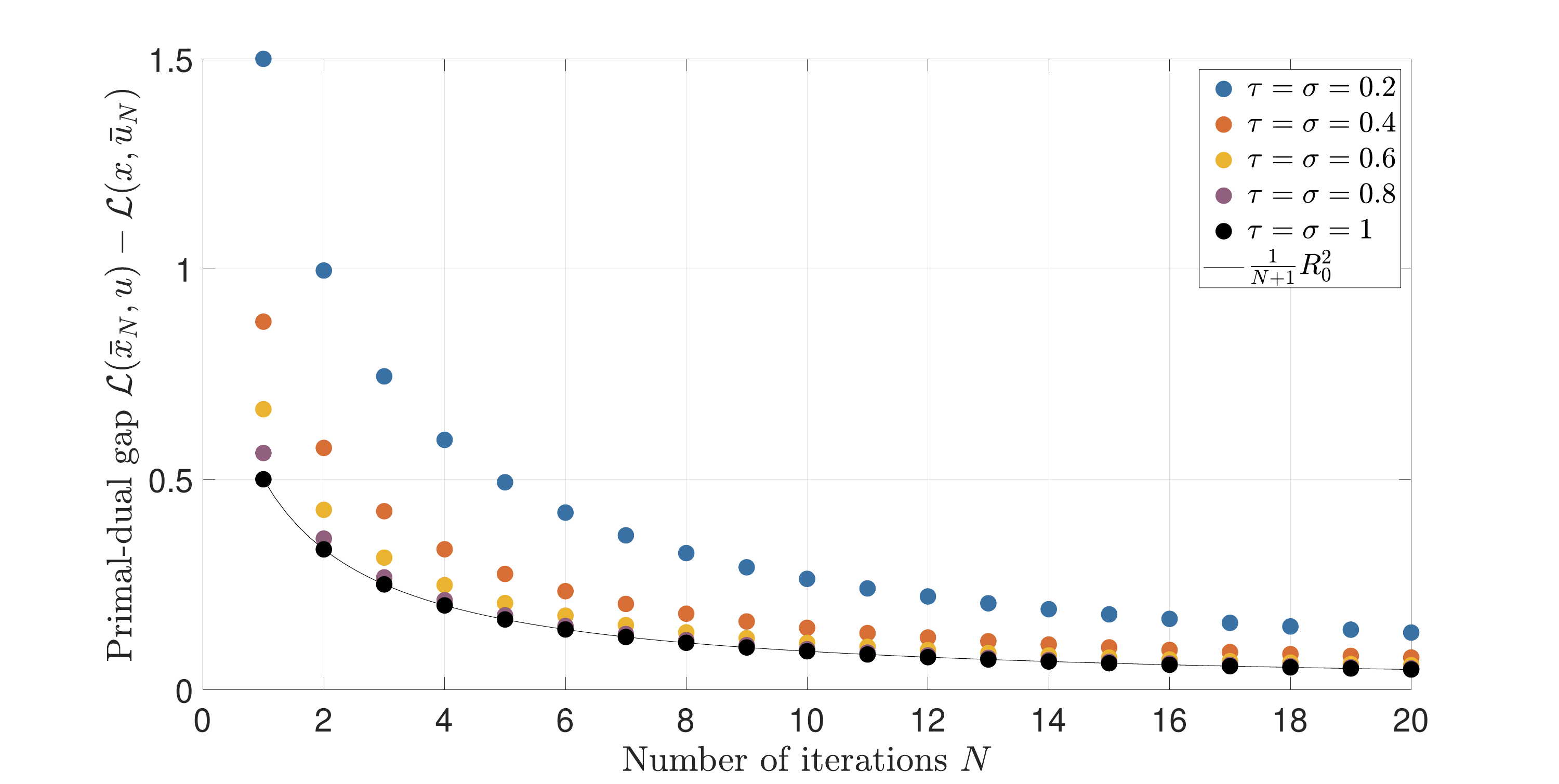}
    }
    \setlength{\abovecaptionskip}{-1pt}
     \setlength{\belowcaptionskip}{-1pt}
    \caption{Worst-case performance obtained by our extension of PEP for $N$ iterations of the Chambolle-Pock algorithm \eqref{eq:(CP)} with different step sizes $\tau =\sigma$ on the problem $\min_x F(x)$ where $F = f + g \circ M$, $f$ and $g$ are convex proximable and $M$ is such that $0\leq ||M||\leq 1$. The performance criterion is the primal-dual gap $\mathcal{L}(\bar{x}_N,u) - \mathcal{L}(x,\bar{u}_N)$ and the initial distance is $R_0^2 = \norm{x-x_0}^2 + \norm{u-u_0}^2$.
    }
    \label{fig:CP16_no_3_term}
\end{figure}

\noindent In this case, our numerical results indicate that the primal-dual gap depends on the step sizes $\tau$ and $\sigma$. Fig. \ref{fig:CP16_no_3_term} shows the worst-case performance returned by PEP with this initial condition for different values of $\tau =\sigma$. When $\tau=\sigma=1$, the numerical results exactly match the curve $\frac{1}{N+1} R_0^2$ for all tested values of $N$, which can be compared to the rate $\frac{1}{N}R_0^2$ proved in \cite{chambolle2016ergodic}. We could also identify the rate for all tested values of $\sigma = \tau \leq 1$ when $N=1$, for which the numerical results match $\frac{\tau+1}{4 \tau}R_0^2$. While we were not able to identify the exact rate for general $N$, $\tau$, and $\sigma$, the expression $\frac{\tau+1}{2\tau(N+1)}R_0^2$ appears to be a valid upper bound when $\tau = \sigma$ (and becomes exact when $\tau=\sigma=1$ or $N=1$ as mentioned above) according to our numerical observations. \vspace{-0.3cm}

\subsubsection{Convergence results on the class of Lipschitz convex functions}
The PEP framework allows us to identify performance guarantees in the setting of our choice, i.e. we can study any performance criterion under any initial condition. 

Theorems \ref{th:CP11} and \ref{th:CP16} characterized the primal-dual gap without any boundedness assumptions on the functions $f$ and $g$. To conclude this section, we propose an alternative analysis, where we consider the primal accuracy, i.e. the value of the (primal) objective function at some iterate. For the corresponding rates to be finite, we must consider a bounded class of functions. We choose the class of Lipschitz convex functions, i.e.\@ convex functions with bounded subgradient. Note that we could have used any other type of bounded classes, e.g. $L$-smooth functions.

We use $||x_0 - x^*||^2 \leq R_x^2$ and $\norm{u_0-u^*}^2 \leq R_u^2$ as a pair of initial conditions. We have to fix the values of all parameters $\tau$, $\sigma$, $N$, $L_M$, $R_x$ and $R_u$, and we consider not necessarily symmetric matrices $M$. 
To show the flexibility of our approach, in addition to computing the rate for the (standard) average iterate, we also perform the analysis for the last iterate, for the best iterate, and for two other averages of iterates, namely the average of the last $\frac{N}{2}$ iterations, and an average using linearly increasing weights (weight $i$ for iteration $i$ before normalization).

\begin{figure}[H]
    \centering
    \figureH{
    \includegraphics[width=1\textwidth]{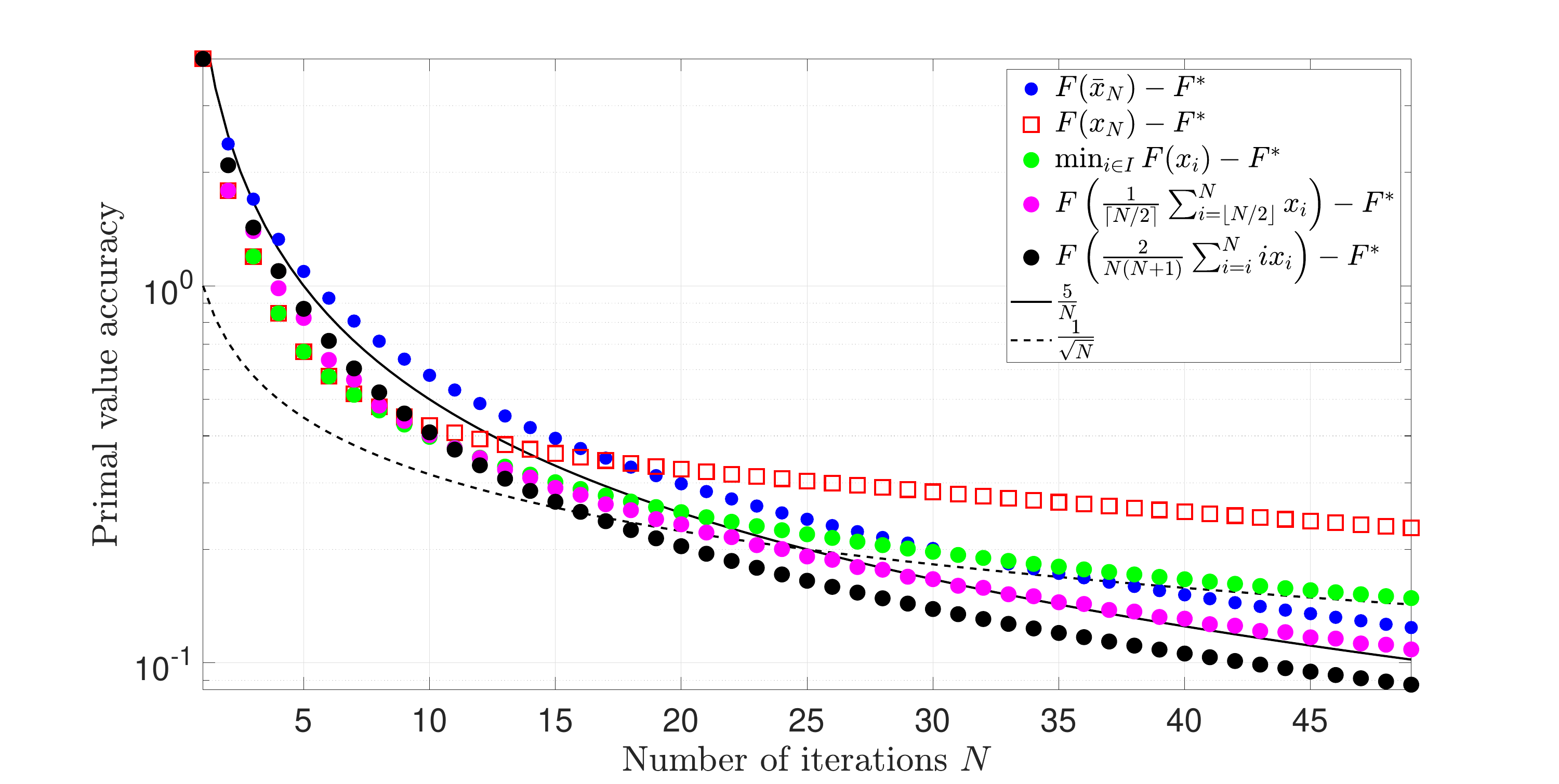}
    }
    \setlength{\abovecaptionskip}{-0pt}
     \setlength{\belowcaptionskip}{-0pt}
    \caption{Worst-case performance obtained by our extension of PEP for $N$ iterations of Chambolle-Pock algorithm \eqref{eq:(CP)} with step size parameters $\tau =\sigma = 1$ on the problem $\min_x F(x)$ where $F = f + g \circ M$, $f$ and $g$ are 1-Lipschitz convex proximable, and $M$ is such that $0\leq ||M||\leq 1$. The performance criterion is the objective function accuracy of the average (blue dots), last (red squares), best (green dots), last $\frac{N}{2}$ (magenta dots), and weighted sum (black dots) of iterates. Curves $\frac 5N$ (solid black line) and $\frac{1}{\sqrt{N}}$ (solid black dashed line) are also represented for comparison purposes.
    }
    \label{fig:CP}
\end{figure}

Fig. \ref{fig:CP} displays the worst-case performance in the above five cases as computed by PEP when minimizing $F = f + g \circ \mathcal{M}$ with \eqref{eq:(CP)} for a number of iterations ranging from $N=1$ to $N=50$. We bounded the primal and dual initial distances $||x_0 - x^*||^2 \leq R_x = 1$ and $||u_0 - u^*||^2 \leq R_u = 1$ and fixed $\tau = \sigma =1$. The primal accuracy of the average iterate (blue dots) seems to roughly follow the $\frac{5}{N}$ curve, whereas the last (red squares) and best (green dots) iterates appear closer to the $\frac{1}{\sqrt{N}}$ curve in this example. Moreover, returning the average of the last $\frac{N}{2}$ iterations (magenta dots) or the proposed weighted sum of the iterates (black dots) leads to even better performances of the method.
\vspace{-0.4cm}

\noindent Numerical guarantees such as those depicted in Fig. \ref{fig:CP} can also be easily obtained for other performance criteria (e.g. primal-dual gap, dual value accuracy), initial distance conditions, and function classes (e.g. symmetric linear operator with lower bounded eigenvalues).  
These could be of great help in the analysis and exploration of the algorithm performance and the identification of interesting phenomena. Moreover, any bound obtained in such manner will be exact, i.e.\@ unimprovable over the considered function class. By contrast, analytical results typically found in the literature may be subject to nontrivial modification or even have to be re-developed when changing the framework of evaluation. For example, the convergence of a weighted sum of the iterates as done in Fig. \ref{fig:CP} would in all likelihood be quite difficult to analyze with standard techniques.

\section{Conclusion}\label{sect:conclusion}
Our main contribution is twofold. First, we obtained interpolation conditions for {classes of linear operators $\LO{L}$, $\LOS{\mu}{L}$, and $\LOT{L}$} with bounds on the eigenvalues or singular values spectrum. 
With these, the powerful and well-developed framework of Performance Estimation Problems can be extended to analyze first-order optimization methods applied to problems involving linear operators.

Our extension allows us to analyze any fixed-step first-order method with any of the standard performance criteria. As a first example, we analyze the worst-case behavior of the gradient method applied to the problem $\min_x g(Mx)$ where $g$ is smooth and convex. We compute exact performance guarantees and conjecture a closed-form expression for the convergence rate. We also analyze the more sophisticated Chambolle-Pock algorithm, tighten the existing primal-dual gap convergence rate, and obtain new exact, numerical convergence guarantees when the objective components are convex and Lipschitz.

Throughout this work, we have used the PEP framework to obtain different types of performance guarantees, ranging from purely analytical to purely numerical, depending on the intrinsic complexity of the desired result and our abilities. More precisely, we encountered the following, distinct situations \begin{enumerate}
    \item  design using PEP-identified coefficients of a new, independent mathematical proof of a performance guarantee (see the proof for one step of the gradient method in Section \ref{sect:perf_gMx}) \item identification of the analytical expression of an exact convergence rate, later proved using an improvement of the  classic argument (see Theorem \ref{th:CP16_improved} on the primal-dual gap of \eqref{eq:(CP)})
    \item identification of new analytical expression for some exact performance guarantees, from which a theoretical explanation is conjectured  (see Conjecture \ref{conj:perf_gMx} on the gradient method)
    \item identification of a new analytical expression for a performance guarantee (see Fig. \ref{fig:CP16_no_3_term} for the case $\sigma =  \tau =1$ on the primal-dual gap of \eqref{eq:(CP)} with another initial distance)
    \item purely numerical performance guarantee (see Fig. \ref{fig:CP} for \eqref{eq:(CP)} on Lipschitz convex functions)
    \end{enumerate}

For future research, we hope to exploit our new approach to improve our understanding of the large variety of methods tackling problems involving linear operators and compare them in terms of worst-case performance.

\section*{Acknowledgments}
The authors would like to thank Prof. Paul Van Dooren for pointing us the existence of some linear algebra results and Anne Rubbens, Moslem Zamani, and Teodor Rotaru for fruitful discussions.

\nocite{nesterov1983method}

\bibliographystyle{apalike} 
\bibliography{sections/references.bib}

\appendix
\section{Properties of matrices}
We review some results of linear algebra used in our proofs. We use $\eye{d}$ and $\zeros{m}{n}$ for identity and zero matrices (dimension may be omitted) and $||M||=\sigma_{\max}(M)$ for the norm of $M$.
\begin{proposition}[\cite{boyd2004convex}; \cite{gallier2010notes}, Theorem 4.3] \label{prop:schur}
Let $G = \vmat{A}{B}{B^T}{C}$ symmetric. We have the three following equivalences
\begin{equation}
    G \succeq 0 \Leftrightarrow 
    \begin{cases}
    C & \succeq 0, \\
    A - B\pinv{C} B^T & \succeq 0, \\
    (I-C\pinv{C} )B^T & = 0,
    \end{cases}\Leftrightarrow 
    \begin{cases}
    A & \succeq 0, \\
    C - B^T\pinv{A} B & \succeq 0, \\
    (I-A\pinv{A} )B & = 0.
    \end{cases}
\end{equation}
\end{proposition}

\begin{proposition}\label{prop:AAB_B}
Let a symmetric matrix $A\succeq 0$. We have $A\pinv{A} = A^\half (\pinv{A})^\half$.
\begin{proof}
We have $A\pinv{A} = A^\half A^\half \pinv{(A^\half)} \pinv{(A^\half)} = A^\half \pinv{(A^\half)} A^\half \pinv{(A^\half)} = A^\half \pinv{(A^\half)}$ where $A^\half \pinv{(A^\half)}  = \pinv{(A^\half)}  A^\half$ holds by definitions of pseudo inverse and square root matrix. 
\end{proof}
\end{proposition}

\begin{proposition}\label{prop:XA_XA}
Let two matrices $C$ and $X$. We have $XC=0 \Leftrightarrow XCC^T = 0$.
\begin{proof}
    If $XCC^T= 0$, then, $XCC^T (\pinv{C})^T =  XC\pinv{C} C = X C = 0$ and if $XC = 0$, then, $XCC^T =0$. 
\end{proof}
    
\end{proposition}

\begin{proposition}\label{prop:null_loewner}
Let two symmetric matrices $A\succeq 0$ and $C\succeq 0$. We have\vspace{-0.2cm}
\begin{equation}
    \exists \alpha>0~:~C \preceq \alpha A \Leftrightarrow A\pinv{A}C=C.
\end{equation}

\vspace{-0.25cm}

\begin{proof}
\proofH{
By application of Proposition \ref{prop:schur}, we have \vspace{-0.2cm}
\begin{align}
    \begin{split}
        \exists \alpha>0~:~
        \begin{cases}
            \alpha I \succ 0, \\
            C \preceq \alpha A,
        \end{cases}
        & \Leftrightarrow 
        \exists \alpha>0~:~\vmat{A}{C^\half}{C^\half}{\alpha I}\succeq 0 \\
        & \Leftrightarrow
        \exists \alpha>0~:~
        \begin{cases}
            A \succeq 0, \\
            \alpha I \succeq C^\half \pinv{A} C^\half, \\
            (I-A\pinv{A})C^\half = 0,
        \end{cases} \\
        & \Leftrightarrow
        (I-A\pinv{A})C = 0,
    \end{split}
\end{align}

\vspace{-0.4cm}

\noindent since $A\succeq 0$, $(I-A\pinv{A})C^\half = 0 \Leftrightarrow (I-A\pinv{A})C = 0$ by Proposition \ref{prop:XA_XA}, and that we can always find an $\alpha$ (sufficiently large) such that $\alpha I \succeq C^\half \pinv{A} C^\half$.  
}
\end{proof}
\end{proposition}

The following result allows to extend a block matrix without increasing its maximal singular value. It will be used to prove Lemma \ref{lem:MR_nonsym} in Appendix \ref{sect:appendix_proof1}.
\begin{theorem}[\cite{davis1982norm}, Theorems 1.1 and 1.2]\label{lem:extension}
Let $M_i$ be conformable matrices. \\
$\exists W$ s.t.  $\norm{\vmat{M_1}{M_2}{M_3}{W}} \leq L$ if, and only if, $\norm{\vrow{M_1}{M_2}} \leq L$ and $\norm{\vcol{M_1}{M_3}} \leq L$.

\noindent Moreover, if $M_1$ is symmetric (resp. skew-symmetric) and $M_2 = M_3^T$ (resp. $M_2 = -M_3^T$), then, there is a $W$ symmetric (resp. skew-symmetric).
\end{theorem}
Note that extending a matrix while maintaining its largest singular value under a given bound $L$ is not trivial in general and $W=0$ does not always work. For example, $\norm{\vrow{1}{1}} = \sqrt{2}$, $\norm{\vcol{1}{1}} = \sqrt{2}$, $\norm{\vmat{1}{1}{1}{0}} = \frac{1+\sqrt{5}}{2} \approx 1.618 > \sqrt{2}$ and $\norm{\vmat{1}{1}{1}{-1}} = \sqrt{2}$, therefore, $W=-1$ is a solution and $W=0$ is not.

The following result states that vector sets leading to the same Gram matrix are equal up to a rotation when the vectors of the two sets have the same dimension. It will be used to prove Lemma \ref{lem:rotation_XYUV} in Appendix \ref{sect:appendix_proof2}.
\begin{theorem}[\cite{horn2012matrix}, Theorem 7.3.11]\label{lem:rotation}
$A$ and $B \in \mathbb{R}^{d \times N}$ build the same Gram matrix, i.e. $A^T A = B^T B$, if, and only if,
\begin{equation}
    \exists V \in \mathbb{R}^{d \times d} \text{ unitary } : B = VA.
\end{equation}
\end{theorem}

\section{Proof of Lemma \ref{lem:MR_nonsym}} \label{sect:appendix_proof1}

\begin{proof}
\proofH{

Let $G = \vmat{A_1}{B_1}{B_1^T}{C_1}$ and $H=\vmat{A_2}{B_2}{B_2^T}{C_2}$ symmetric and positive semidefinite matrices satisfying \vspace{-0.2cm}
\begin{equation}\label{eq:cond_ABC}
    \begin{cases}
        B_1 = B_2, \\
        A_2 \preceq L^2 A_1, \\
        C_1 \preceq L^2 C_2.
    \end{cases}
\end{equation}
In the remaining of the proof, we equivalently use $B_1$ or $B_2$. We note $S_1 = C_1 - B_1^T \pinv{A_1}B_1$ and $S_2 = A_2-B_2\pinv{C_2}B_2^T $ where $A_1,C_1,A_2,C_2,S_1,S_2 \succeq 0$ by Proposition \ref{prop:schur} and $G,H\succeq 0$. Moreover, Proposition \ref{prop:schur} and \ref{prop:AAB_B} with $G,H \succeq 0$ provide
\begin{align}\label{eq:eq_for_proof_lemma35}
    \begin{split}
    A_1^\half (\pinv{A_1})^\half B_1 = A_1 \pinv{A_1} B_1 = B_1, \\
    B_2 (\pinv{C_2})^\half C_2^\half = B_2 \pinv{C_2} C_2 = B_2,
    \end{split}
\end{align}
and Propositions \ref{prop:AAB_B} and \ref{prop:null_loewner} with \eqref{eq:cond_ABC} and \eqref{eq:eq_for_proof_lemma35} provide
\begin{align}\label{eq:cond_SABC}
\begin{split}
    S_2^\half A_1^\half (\pinv{A_1})^\half = S_2^\half, \\
    S_1^\half C_2^\half (\pinv{C_2})^\half = S_1^\half.
    \end{split}
\end{align}

\vspace{-0.4cm}

Firstly, we show that \vspace{-0.2cm}
\begin{equation}
    (X_R,Y_R,U_R,V_R) = \left( \vcol{A_1^\half}{\zeros{N_1}{N_1}}, \vcol{(\pinv{C_2})^\half B_2^T}{S_2^\half},\vcol{C_2^\half}{\zeros{N_1}{N_2}}, \vcol{(\pinv{A_1})^\half B_1}{S_1^\half}  \right)
\end{equation} 

\noindent is a factorization of $G$ and $H$. Indeed, we have
\vspace{-0.2cm}
\begin{align}
    X_R^T X_R & = \vrow{A_1^\half}{\zeros{N_1}{N_2}}\vcol{A_1^\half}{\zeros{N_2}{N_1}} = A_1, \\
    U_R^T U_R & = \vrow{C_2^\half}{\zeros{N_2}{N_1}}\vcol{C_2^\half}{\zeros{N_1}{N_2}} = C_2, \\
    Y_R^T Y_R & = \vrow{B_2(\pinv{C_2})^\half}{S_2^\half}\vcol{(\pinv{C_2})^\half B_2^T}{S_2^\half} = B_2 (\pinv{C_2})^\half (\pinv{C_2})^\half B_2^T + S_2 = A_2, \\
    V_R^T V_R & = \vrow{B_1^T(\pinv{A_1})^\half}{S_1^\half}\vcol{(\pinv{A_1})^\half B_1}{S_1^\half} = B_1^T (\pinv{A_1})^\half (\pinv{A_1})^\half B_1 + S_1 = C_1, \\
    X_R^T V_R & = \vrow{A_1^\half}{\zeros{N_1}{N_2}}\vcol{(\pinv{A_1})^\half B_1}{S_1^\half} = \overbrace{A_1^\half (\pinv{A_1})^\half B_1}^{=B_1 \text{ by } \eqref{eq:eq_for_proof_lemma35}} = B_1, \\
    Y_R^T U_R & = \vrow{B_2(\pinv{C_2})^\half}{S_2^\half}\vcol{C_2^\half}{\zeros{N_1}{N_2}} = \overbrace{B_2(\pinv{C_2})^\half C_2^\half}^{=B_2 \text{ by } \eqref{eq:eq_for_proof_lemma35}} = B_2,
\end{align}
and therefore
\begin{align}
\begin{split}
    \vrow{X_R}{V_R}^T\vrow{X_R}{V_R} & = \vmat{X_R^TX_R}{X_R^TV_R}{V_R^TX_R}{V_R^TV_R}  = \vmat{A_1}{B_1}{B_1^T}{C_1} = G, \\
    \vrow{Y_R}{U_R}^T\vrow{Y_R}{U_R} & = \vmat{Y_R^TY_R}{Y_R^TU_R}{U_R^T Y_R}{U_R^T U_R}  = \vmat{A_2}{B_2}{B_2^T}{C_2} = H.
\end{split}
\end{align}

Secondly, we show that $(X_R,Y_R,U_R,V_R)$ is $\LO{L}$-interpolable by providing a linear operator $M_R$ such that $Y_R = M_R X_R$, $V_R = M_R^T U_R$ and $||M_R|| \leq L$. The expression of $M_R$ is
\begin{equation}\label{eq:MR2}
    M_R = \vmat{(\pinv{C_2})^\half B_1^T (\pinv{A_1} )^\half}{(\pinv{C_2} )^\half S_1^\half}{ S_2^\half (\pinv{A_1} )^\half}{W}
\end{equation}
where $W$ is a matrix specified later. Indeed, we have
\begin{align}
\begin{split}
    M_R X_R & = \vmat{(\pinv{C_2})^\half B_1^T (\pinv{A_1} )^\half}{(\pinv{C_2} )^\half S_1^\half}{ S_2^\half (\pinv{A_1} )^\half}{W} \vcol{A_1^\half}{\zeros{N_2}{N_1}} \\
    & = \vcol{(\pinv{C_2})^\half B_1^T (\pinv{A_1} )^\half A_1^\half}{ \underbrace{S_2^\half (\pinv{A_1} )^\half A_1^\half}_{=S_2^\half \text{ by } \eqref{eq:cond_SABC}}} = \vcol{(\pinv{C_2})^\half B_1^T}{ S_2^\half} = Y_R, 
\end{split}
\end{align}
\begin{align}
\begin{split}
    M_R^T U_R & = \vmat{(\pinv{A_1})^\half B_2^T (\pinv{C_2} )^\half}{(\pinv{A_1} )^\half S_2^\half}{ S_1^\half (\pinv{C_2} )^\half}{W^T} \vcol{C_2^\half}{\zeros{N_1}{N_2}} \\
    & = \vcol{(\pinv{A_1})^\half B_2^T (\pinv{C_2} )^\half C_2^\half}{\underbrace{S_1^\half (\pinv{C_2} )^\half C_2^\half}_{=S_1^\half \text{ by } \eqref{eq:cond_SABC}}} = \vcol{(\pinv{A_1})^\half B_2^T}{S_1^\half} = V_R.
\end{split}
\end{align}

It remains to show that the proposed $M_R$ has singular values bounded by $L$ for a suited choice of $W$. Thanks to Theorem \ref{lem:extension}, we must just show that the matrices 
\begin{align}
\begin{split}
    M_R^{(\text{up})} & = \vrow{(\pinv{C_2} )^\half B_1 (\pinv{A_1} )^\half}{(\pinv{C_2} )^\half S_1^\half}, \\
    M_R^{(\text{left})} & = \vcol{(\pinv{C_2} )^\half B_2^T (\pinv{A_1} )^\half}{S_2^\half (\pinv{A_1} )^\half},
\end{split}
\end{align}

have singular values lower than $L$, or equivalently, that the products $M_R^{(\text{up})} {M_R^{(\text{up})}}^T$ and ${M_R^{(\text{left})}}^T M_R^{(\text{left})}$ have eigenvalues lower than $L^2$, i.e.,\vspace{-0.2cm}
\begin{align}
\begin{split}
    M_R^{(\text{up})} {M_R^{(\text{up})}}^T & = (\pinv{C_2} )^\half B_1 \pinv{A_1}  B_1^T (\pinv{C_2} )^\half + (\pinv{C_2} )^\half S_1 (\pinv{C_2} )^\half\\
    & = (\pinv{C_2} )^\half C_1 (\pinv{C_2} )^\half \preceq L^2 I 
\end{split}
\end{align}
\vspace{-0.2cm} \noindent
and \vspace{-0.2cm}
\begin{align}
\begin{split}
     {M_R^{(\text{left})}}^T M_R^{(\text{left})} & = (\pinv{A_1} )^\half B_2 \pinv{C_2}  B_2^T (\pinv{A_1} )^\half + (\pinv{A_1} )^\half S_2 (\pinv{A_1} )^\half \\
     & = (\pinv{A_1} )^\half A_2 (\pinv{A_1} )^\half \preceq L^2 I 
 \end{split}
\end{align}

which are both true since $C_1 \preceq L^2 C_2$ implies that \vspace{-0.2cm}
\begin{align*}
        (\pinv{C_2})^\half C_1 (\pinv{C_2})^\half & \preceq  L^2 (\pinv{C_2})^\half C_2^\half C_2^\half (\pinv{C_2})^\half  = L^2 C_2 \pinv{C_2} C_2 \pinv{C_2}  = L^2 C_2 \pinv{C_2}  \preceq L^2I
\end{align*}

by definition of the pseudo inverse and Proposition \ref{prop:AAB_B} (same result on $A_2\preceq L^2 A_1$).

Finally, we observe on expression \eqref{eq:MR2} of $M_R$ that, if $A_1=C_2$ and $A_2=C_1$, and $B=B^T$ (resp. $B=-B^T$), then thanks to Theorem \ref{lem:extension}, we can choose $W$ symmetric (resp. skew-symmetric) such that $M_R$ is symmetric (resp. skew-symmetric). Moreover, if $A_1=C_2$, then, $U_R = X_R$ and $V_R = Y_R$ (resp. $V_R = -Y_R$). Note that in the skew-symmetric case, we have to add a negative sign on one of the two off-diagonal blocks of $M_R$ in \eqref{eq:MR2}.
}
\end{proof}

\section{Proof of Lemma \ref{lem:rotation_XYUV}} \label{sect:appendix_proof2}

\begin{proof}
    \proofH{
Adding zeros to $\vrow{X}{V}$ or $\vrow{X_R}{V_R}$ such that they have the same number of rows, i.e. $d_n = \max\{n,n_R\}$, will allow us to use Theorem \ref{lem:rotation}. We use $E_{n,d_n} = \vcol{\eye{n}}{0_{(d_n-n),n}}$ and $E_{n_R,d_n} = \vcol{\eye{n_R}}{0_{(d_n-n_R),n_R}}$ to add the zeros and now work with $E_{n,d_n} \vrow{X}{V}$ and $E_{n_R,d_n} \vrow{X_R}{V_R}$.
Moreover, adding the zeros preserves the Gram matrices, indeed, if $\vrow{X}{V}^T \vrow{X}{V} = \vrow{X_R}{V_R}^T \vrow{X_R}{V_R}$, then \vspace{-0.15cm}
    \begin{align}
        \vrow{X}{V}^T E_{n,d_n}^T E_{n,d_n} \vrow{X}{V} & = \vrow{X_R}{V_R}^T E_{n_R,d_n}^T E_{n_R,d_n} \vrow{X_R}{V_R}.
    \end{align}
    \vspace{-0.6cm}
    
\noindent Therefore, Theorem \ref{lem:rotation} applies and yields \vspace{-0.1cm}
    \begin{align}
        E_{n_R,d_n} \vrow{X_R}{V_R} & = V_G E_{n,d_n} \vrow{X}{V}
    \end{align}
    \vspace{-0.6cm}
    
\noindent and with the same reasoning on $\vrow{Y}{U}$, $\vrow{Y_R}{U_R}$, and $d_m = \max\{m,m_R\}$, it yields \vspace{-0.4cm}
    \begin{align}
        E_{m_R,d_m} \vrow{Y_R}{U_R} & = V_H E_{m,d_m} \vrow{Y}{U}
    \end{align}
    \vspace{-0.6cm}
    
\noindent for some $V_G$ and $V_H$ unitary.
Now, since  $(X_R,Y_R,U_R,V_R)$ is {$\LO{L}$}-interpolable, there exists a $M_R$ such that $Y_R = M_R X_R$, $V_R = M^T U_R$, and $||M_R||\leq L$. Therefore,
\begin{align}
    Y_R  = M_R X_R  & \Rightarrow \overbrace{E_{m_R,d_m}Y_R}^{V_H E_{m,d_m} Y} =  E_{m_R,d_m} M_R E_{n_R,d_n}^T \overbrace{E_{n_R,d_n} X_R}^{V_G E_{n,d_n} X} \\
                    & \Rightarrow Y  = \overbrace{E_{m,d_m}^T V_H^T E_{m_R,d_m} M_R E_{n_R,d_n}^T V_G E_{n,d_n}}^{M} X 
\end{align}
and
\begin{align}
    V_R  = M_R^T U_R  & \Rightarrow \overbrace{E_{n_R,d_n}V_R}^{V_G E_{n,d_n} V} =  E_{n_R,d_n} M_R^T E_{m_R,d_m}^T \overbrace{E_{m_R,d_m} U_R}^{V_H E_{m,d_m} X} \\
                    & \Rightarrow V  = \overbrace{E_{n,d_n}^T V_G^T E_{n_R,d_n} M_R^T E_{m_R,d_m}^T V_H E_{m,d_m}}^{M^T} U 
\end{align}
We have $\norm{M}\leq L$ since unitary transformations $V_H^T$ and $V_G$ preserve the maximal singular value and that $E_{i,j}$ can only add zeros singular values.\\
Finally, when $U=X$, $V=Y$ (resp. $V=-Y$), $U_R=X_R$ and $V_R=Y_R$ (resp. $V_R = -Y_R$), we have $V_G=V_H$. Therefore, $M$ is obtained as a unitary transformation of $M_R$, in other words, if $M_R$ was symmetric (resp. skew-symmetric), then, $M$ remains symmetric (resp. skew-symmetric).} 
\end{proof}

\section{Proof of the equality in \eqref{eq:proof_1_step}}\label{app:proof_1_step}
The equality in \eqref{eq:proof_1_step} holds when $\gamma = \frac{1-\mu_g}{2-\mu_g}$, $\tau = \frac{\mu_g}{2(\mu_g-1+(1-\mu_g)^{-2})}$, $x_1 = x_0 - M^T \nabla g(Mx_0)$, and
\begin{equation*}
    \mu_g \leq \frac16 \left(7 - \frac{7}{\sqrt[3]{44 - 3 \sqrt{177}}} - \sqrt[3]{44 - 3 \sqrt{177}} \right) \approx 0.17,
\end{equation*}
\begin{align*}
    P & = \begin{pmatrix} Mx_0 & Mx_1 & Mx^* & \nabla g(M x_0) & \nabla g(M x_1) & \nabla g(M x^*)\end{pmatrix},
\end{align*}
\begin{equation*}
    a_{\mu_g}  = \begin{pmatrix} a_1 \\ -(1+\mu_g)a_2 \\ a_3 \\ -(1+\mu_g)a_2 \\ a_2 \\ \mu_g a_2 \end{pmatrix},
    b_{\mu_g}  = 
    \begin{pmatrix}
        0 \\b_1 \\b_2 \\b_3 \\ b_4 \\ b_5 
    \end{pmatrix}=
    \begin{pmatrix} 0 \\
                    -\sqrt{\frac{\mu_g}{2(1-\mu_g)} + \frac{1}{2-\mu_g} + (1+\mu_g)^2 a_2^2} \\
                    -b_1 \\
                    \frac{2 (1+\mu_g)^2 (2-\mu_g) a_2^2 - (1+\mu_g)}{-2(2-\mu_g) b_1} \\
                    \frac{-2 (1-\mu_g^2) a_2^2 + 1}{-2 (1-\mu_g) b_1} \\
                    \frac{ -2\mu_g(1-\mu_g^2) (2-\mu_g) a_2^2 +(1-\mu_g)^2 + \mu_g }{2 (1-\mu_g) (2-\mu_g) b_1}         
    \end{pmatrix},
\end{equation*}
\begin{equation*}
    c_{\mu_g}  = 
    \begin{pmatrix} 0 \\ 0 \\ 0 \\ \sqrt{ \frac{1}{2-\mu_g} - (1+\mu_g)^2 a_2^2 - b_2^2} \\ - \sqrt{\frac{1}{2(1-\mu_g)} -a_2^2 -b_3^2} \\  \sqrt{\frac{1}{2(1-\mu_g)} - (\mu_g a_2)^2 - b_4^2}\end{pmatrix} 
\end{equation*}
where $ a_1 = -\sqrt{\frac{\mu_g}{2-\mu_g} + \tau}$, $a_2 = \frac{1}{2(2-\mu_g)a_1}$, $a_3 = \sqrt{\tau - \frac{1}{2-\mu_g} + (1+\mu_g)^2 a_2^2 }$. Indeed, given these definitions, elementary algebraic calculations allow to prove the equality.

\section{Proof of Theorem \ref{th:Beck}} \label{sect:Beck}
The proof comes from a lecture by Prof. Amir Beck \cite{beck2022}.
\begin{proof}
    Since the assumptions of Theorem \ref{th:CP16} hold, it applies $\forall x,y$, namely,
    \begin{equation*}
    \mathcal{L}(\bar{x}_N,u) - \mathcal{L}(x,\bar{u}_N) \leq  \frac{1}{2N} \Par{\frac{\norm{x-x_0}^2}{\tau} + \frac{\norm{u-u_0}^2}{\sigma} - 2 (u-u_0)^TM(x-x_0)}
    \end{equation*}
    where $\mathcal{L}(x,u) = x^T M^T u - g^*(u) + f(x)$. Choosing $u = \tilde{u}_N \in \partial g(M\bar{x}_N)$, $x = x^*$ and using Conjugate Subgradient Theorem (CST) and Fenchel's Inequality (FI) yield
    \begin{align*}
        \mathcal{L}(\bar{x}_N,\tilde{u}_N) & = \underbrace{\bar{x}_N^T M^T \tilde{u}_N - g^*(\tilde{u}_N)}_{=g(M\bar{x}_N)\text{ by (CST)}} +f(\bar{x}_N) =  g(M\bar{x}_N) + f(\bar{x}_N), \\
        \mathcal{L}(x^*,\bar{u}_N) & = \underbrace{{x^*}^T M^T \bar{u}_N -g^*(\bar{u}_N)}_{\leq g(Mx^*)\text{ by (FI)}} + f(x^*) \leq g(Mx^*) + f(x^*),
    \end{align*}
    and therefore 
    \begin{align*}
        F(x_N) - F(x^*) & \leq \mathcal{L}(\bar{x}_N,\tilde{u}_N) -\mathcal{L}(x^*,\bar{u}_N) \\
        & \leq \frac{1}{2N} \Par{\frac{\norm{x^*-x_0}^2}{\tau} + \frac{\norm{\bar{u}_N-u_0}^2}{\sigma} - 2 (\bar{u}_N-u_0)^TM(x^*-x_0)}
    \end{align*}
    where $F(x) = f(x) + g(Mx)$.
\end{proof}

\end{document}